\documentclass{amsart}

\newtheorem{theorem}{Theorem}[section]%[subsection]
\newtheorem{lemma}[theorem]{Lemma}
\newtheorem{proposition}[theorem]{Proposition}
\newtheorem{corollary}[theorem]{Corollary}
\theoremstyle{definition}
\newtheorem{definition}[theorem]{Definition}
\newtheorem{example}[theorem]{Example}

\newtheorem{definitionandtheorem}[theorem]{Definition and Theorem}

\theoremstyle{remark}
\newtheorem{remark}[theorem]{Remark}

\newcommand{\Tr}{\mbox{Tr}} 
\newcommand{\tr}{\mbox{tr}}
\newcommand{\id}{\mbox{id}}

\newcommand{\End}{\mbox{End}} 
\newcommand{\Rep}{\mbox{Rep}}
\newcommand{\URep}{\mbox{URep}}

\newcommand{\Hom}{\mbox{Hom}}

\renewcommand{\span}{\mbox{span}} 
    
\newcommand{\Ad}{\mbox{Ad}\,}
\newcommand{\R}{{\mathcal R}}
\newcommand{\bR}{{\bar \R}}

\newcommand{\C}{{\mathcal C}}

\newcommand{\E}{{\varepsilon}}

\newcommand{\eps}{\varepsilon} 
\newcommand{\op}{op}
\newcommand{\cop}{cop}

\newcommand{\actr}{\rightharpoonup} 
\newcommand{\actl}{\leftharpoonup}

\newcommand{\la}{\langle\,} 
\newcommand{\ra}{\,\rangle}

\newcommand{\rtimes}{{>\!\!\!\triangleleft}} 
\newcommand{\ltimes}{{\triangleright\!\!\!<}}

\newcommand{\own}{\ni} 
  
\newcommand{\1}{_{(1)}} 
\newcommand{\2}{_{(2)}} 
\newcommand{\3}{_{(3)}} 
\newcommand{\4}{_{(4)}} 
\newcommand{\5}{_{(5)}}
\newcommand{\6}{_{(6)}}

\newcommand{\I}{^{(1)}} 
\newcommand{\II}{^{(2)}}

\begin{document}

\title{\bf Invariants of Knots and 3-manifolds from Quantum Groupoids}

% Remove or comment out any unused author tags.
% author one information
\author{Dmitri Nikshych}
\address{UCLA, Department of Mathematics, 405 Hilgard Avenue,
Los Angeles, CA 90095-1555, USA}
%\curraddr{}
\email{nikshych@math.ucla.edu}
\thanks{The first author thanks P.~Etingof for useful discussions}

% author two information
\author{Vladimir Turaev}
\address{Institut de Recherche Math\'ematique Avanc\'ee, Universit\'e 
Louis Pasteur,
CNRS, 7, rue Ren\'e Descartes, F-67084 Strasbourg, France}
%\curraddr{}
\email{turaev@math.u-strasbg.fr}
%\thanks{}

% author three information
\author{Leonid Vainerman}
\address{ D\'epartement de
Math\'ematiques, Universit\'e Louis Pasteur, 7, rue Ren\'e Descartes, 
F-67084 Strasbourg, France}
%\curraddr{}
\email{wain@agrosys.kiev.ua, vaynerma@math.u-strasbg.fr}
\thanks{ The third author is grateful to 
l'Universit\'e Louis Pasteur (Strasbourg) 
for the kind hospitality during his work on this article}

%\subjclass{}
\date{June 8, 2000}

\begin{abstract}
We use the categories of representations of finite dimensional quantum
groupoids (weak Hopf algebras) to construct ribbon and modular categories
that give rise to invariants of knots and 3-manifolds.
\end{abstract}

\maketitle
\tableofcontents

\begin{section}
{Introduction}
In \cite{RT2} a general method of constructing invariants of
$3$-manifolds from modular Hopf algebras was introduced. 
After appearance of \cite{RT2}
it became clear that the technique of Hopf algebras can
be replaced by a more general technique of
monoidal categories. An appropriate class of categories 
-- modular categories -- 
was introduced in \cite{T1}. In addition to quantum groups,
such categories also arise from skein categories of tangles
and, as it was observed by A. Ocneanu, from certain bimodule
categories of type II${}_1$
subfactors.

The goal of this paper is to study the representation
categories of {\em quantum groupoids} and to give in this way a new
construction of modular categories.
This extends the construction of
modular categories from modular Hopf algebras and in particular from
quantum groups at roots of unity.

By quantum groupoids, we understand
weak Hopf algebras introduced in \cite{BNSz}, \cite{BSz1}, \cite{Ni}.
These objects generalize Hopf algebras, usual
finite groupoid algebras and their duals (cf. \cite{NV1}).
We use the term ``quantum groupoid'' rather than ``weak Hopf algebra''.

It was shown in \cite{NV2}, \cite{NV3} that 
quantum groupoids and their coideal subalgebras are closely related to
II$_1$-subfactors. Every finite index and finite depth II$_1$-subfactor
gives rise to a pair consisting of a $C^*$-quantum groupoid and its 
left coideal subalgebra, and vice versa.
It was also explained in \cite{NV3}  how to express 
the known subfactor invariants such as bimodule   
categories and principal graphs in terms of
the associated quantum groupoids. 
In particular, the bimodule
categories arising from a finite index and finite depth II$_1$ subfactors 
are equivalent to the unitary representation categories of the corresponding
$C^*$-quantum groupoids
(\cite{NV3}, 5.8).

Thus, it is natural to study categories of representations of quantum 
groupoids and to extend concepts known for Hopf algebras to this setting.
We show that the representation category $\Rep(H)$
of a quantum groupoid $H$ is a monoidal category with duality. We introduce
quasitriangular, ribbon, and modular quantum groupoids for which
$\Rep(H)$ is, respectively, braided, ribbon, and modular.
The notion of factorizability is extended from the Hopf algebra case
and used to construct modular categories. We define
the Drinfeld double $D(H)$ of a quantum groupoid $H$ and show  
that it is a factorizable quasitriangular quantum groupoid. 
For a $C^*$-quantum groupoid $H$, we similarly study the unitary representation
category $\URep(H)$.

It should be mentioned that the category $\URep(H)$ for
a $C^*$-quantum groupoid $H$ was previously introduced by 
G.~B\"ohm and K.~Szlach\'anyi in \cite{BSz2}; they also introduced the 
notion of an $R$-matrix and the Drinfeld double for $C^*$-quantum groupoids, 
see \cite{BSz1}.

Our main theorem (Theorem~\ref{rep of C* is modular}) reads: If $H$ is a
connected 
$C^*$-quantum groupoid, then the category $\URep(D(H))$ of unitary
representations
of $D(H)$ is a 
unitary modular category.

Thus, any finite index and finite depth
II$_1$-subfactor yields a unitary 
modular category as follows: consider the associated connected 
$C^*$-quantum groupoid $H$, then the category
$\URep(D(H))$ is a unitary modular category. 
We conjecture that this construction is equivalent to the one due 
to A.~Ocneanu (see \cite{EK}).
The key role in the proof of the main theorem is played by the following 
Lemma~\ref{factorizable implies modular}:
If $H$ is a connected, ribbon and factorizable quantum groupoid 
with a Haar measure over an algebraically closed field,   
then $\Rep(H)$ is a modular category.

The organization of the paper is clear from the table of contents.
\end{section}
%%%%%%%%%%%%%%%%%%%%%%%%%%%%%%%%%%%%%%%%%%%%%%%%%%%%%%%%%%%%%%%%%%%%%%%%%%%%
%%%%%%%%%%%%%%%     Quantum groupoids             %%%%%%%%%%%%%%%%%%%%%%%%%%
%%%%%%%%%%%%%%%%%%%%%%%%%%%%%%%%%%%%%%%%%%%%%%%%%%%%%%%%%%%%%%%%%%%%%%%%%%%%

\begin{section}
{Quantum groupoids}

In this section we recall basic
properties of quantum groupoids. Most of the material 
presented here can be found in \cite{BNSz} and \cite{NV1}, see also
the survey \cite{NV4}.

Throughout this paper we use Sweedler's notation for
comultiplication, writing $\Delta(b) = b\1 \otimes b\2$.
Let $k$ be a field.

\begin{definition} 
\label{finite quantum groupoid}
A {\em (finite) quantum groupoid} over $k$ is a
finite dimensional $k$-vector space $H$
with the structures of an associative algebra $(H,\,m,\,1)$ 
with multiplication $m:H\otimes_k H\to H$ and unit $1\in H$ and a 
coassociative coalgebra $(H,\,\Delta,\,\eps)$ with comultiplication
$\Delta:H\to H\otimes_k H$ and counit $\eps:H\to k$ such that:
\begin{enumerate}
\item[(i)] The comultiplication $\Delta$ 
is a (not necessarily unit-preserving) homomorphism of algebras such that
\begin{equation}
(\Delta \otimes \id) \Delta(1) =
(\Delta(1)\otimes 1)(1\otimes \Delta(1)) =
(1\otimes \Delta(1))(\Delta(1)\otimes 1),
\label{Delta 1}
\end{equation}
\item[(ii)] The counit is a $k$-linear map
satisfying the identity:
\begin{equation}
\eps(fgh) = \eps(fg\1)\,\eps(g\2h) = \eps(fg\2)\,\eps(g\1h),
\label{eps m}
\end{equation}
for all $f,g,h\in H$. 
\item[(iii)]
There is an algebra and 
coalgebra anti-homomorphism $S: H \to H$, called an {\em antipode},
such that, for all $h\in H$,
\begin{eqnarray}
m(\id \otimes S)\Delta(h) &=&(\eps\otimes\id)(\Delta(1)(h\otimes 1)),
\label{S epst} \\
m(S\otimes \id)\Delta(h) &=& (\id \otimes \eps)((1\otimes h)\Delta(1)).
\label{S epss}
\end{eqnarray}
\end{enumerate} 
\end{definition}
A quantum groupoid is a Hopf algebra if and only if one of the following
conditions holds: (i) the comultiplication 
is unit-preserving or (ii) if and only if the counit is a homomorphism of
algebras.

A {\em morphism} of quantum groupoids is a map between them
which is both an algebra and a coalgebra morphism preserving unit and
counit and commuting with
the antipode. The image of such a morphism is clearly
a quantum groupoid. The tensor product of two quantum groupoids
is defined in an obvious way.

The set of axioms of Definition~\ref{finite quantum groupoid} is self-dual.
This allows to define a natural quantum groupoid
structure on the dual vector space $\widehat H=\Hom_k(H,k)$ by
``reversing the arrows'':
\begin{eqnarray}
& & \la h,\,\phi\psi \ra = \la \Delta(h),\,\phi\otimes\psi \ra, \\
& & \la g\otimes h,\,{\widehat\Delta}(\phi) \ra =  \la gh,\, \phi\ra, \\
& & \la h,\, {\widehat S}(\phi) \ra = \la S(h),\,\phi \ra,
\end{eqnarray}
for all $\phi,\psi \in \widehat H,\, g,h\in H$. The unit $\widehat 1\in
\widehat H$ is $\eps$  and counit $\widehat\eps$ is 
$\phi \mapsto \la\phi,\, 1\ra$.

The linear endomorphisms of $H$ defined by 
\begin{equation}
h\mapsto m(\id \otimes S)\Delta(h), \qquad
h\mapsto m(S \otimes \id)\Delta(h)
\end{equation}
are called the {\em target} and {\em source counital maps} and 
denoted $\eps_t$ and $\eps_s$, respectively. 
From axioms (\ref{S epst}) and (\ref{S epss}),
\begin{equation}
\eps_t(h) = (\eps\otimes\id)(\Delta(1)(h\otimes 1)),\qquad
\eps_s(h) = (\id \otimes \eps)((1\otimes h)\Delta(1)).
\end{equation}
In the Hopf algebra case $\eps_t(h) = \eps_s(h) = \eps(h) 1$. 

We have $S\circ \eps_s = \eps_t\circ S$ and $\eps_s\circ S = S\circ \eps_t$. 
The images of these maps
$\eps_t$ and $\eps_s$
\begin{eqnarray}
H_t &=& \eps_t(H)  =  \{h\in H \mid \Delta(h) = \Delta(1)(h\otimes 1)\}, \\
H_s &=& \eps_s(H)  = \{h\in H \mid \Delta(h) = (1\otimes h)\Delta(1)\}
\end{eqnarray}
are subalgebras of $H$, called the {\em target} (resp.\ {\em source}) 
{\em counital subalgebras}. They play the role of ground
algebras for $H$. They commute with each other and
\begin{equation*}
H_t = \{(\phi\otimes \id)\Delta(1) \mid \phi\in \widehat H \}, \qquad
H_s = \{(\id \otimes \phi)\Delta(1) \mid \phi\in \widehat H \},
\end{equation*}
i.e., $H_t$ (resp.\ $H_s$) is generated by the right (resp.\ left)
tensorands of $\Delta(1)$.  The restriction of $S$
defines an algebra anti-isomorphism between $H_t$ and $H_s$. 
Any non-zero morphism $H \to K$ of  quantum groupoids preserves counital 
subalgebras,  i.e., $H_t \cong K_t$ and $H_s \cong K_s$.

In what follows we will use the Sweedler arrows, 
writing for all $h\in H,\phi\in \widehat H$:
\begin{equation}
h\actr\phi = \phi\1 \la h,\, \phi\2\ra,
\qquad 
\phi\actl h =\la h,\,\phi\1 \ra \phi\2
\end{equation}
for all $h\in H,\phi\in \widehat H$. Then the map    
$z \mapsto (z\actr \eps)$ is an algebra isomorphism
between $H_t$ and ${\widehat H}_s$.
Similarly, the map $y \mapsto (\eps\actl y)$
is an algebra isomorphism between $H_s$ and ${\widehat H}_t$ 
(\cite{BNSz}, 2.6). Thus, the counital subalgebras of $\widehat H$ 
are canonically anti-isomorphic to those of $H$.

A quantum groupoid $H$ is called {\em connected} if
$H_s \cap Z(H) = k$, or, equivalently, $H_t \cap Z(H) = k$, where
$Z(H)$ denotes the center of $H$ (cf. \cite{N}, 3.11, \cite{BNSz}, 2.4).

Let us recall that a $k$-algebra $A$ is {\em separable} \cite{P}
if the multiplication epimorphism $m: A\otimes_k A \to A$
has a right inverse as an $A-A$ bimodule homomorphism.
This is equivalent to the existence of a {\em separability element}
$e\in A\otimes_k A $ such that  $m(e) =1$ and 
$(a\otimes 1)e = e(1\otimes a),\ (1\otimes a)e=e(a\otimes 1)$
for all  $a\in A$.

The counital subalgebras $H_t$ and $H_s$ are separable,
with  separability elements  $e_t = (S \otimes \id)\Delta(1)$ 
and $e_s = (\id \otimes S)\Delta(1)$, respectively.

Observe that the {\em adjoint} actions of $1\in H$ give rise to
non-trivial maps $H\to H$:
\begin{equation}
h \mapsto 1\1hS(1\2) = \Ad_1^l(h), \qquad 
h\mapsto S(1\1) h 1\2 = \Ad_1^r(h), \qquad h\in H.
\end{equation}

\begin{lemma}
\label{expectations}
The map $\Ad^{l}_1$ is a linear projection
from $H$ onto $C_H(H_s)$, the centralizer of $H_s$,
i.e.,\ $(\Ad^{l}_1)^2= \Ad^{l}_1$.
The map $\Ad^{r}_1$ is a linear projection
from $H$ onto $C_H(H_t)$, the centralizer of
$H_t$, i.e.,\ $(\Ad^{r}_1)^2= \Ad^{r}_1$.
\end{lemma}
\begin{proof}
Since $1\1 \otimes S(1\2)$ is a separability element of $H_s$,
$\Ad^{l}_1(h)$ commutes with $H_s$. The assertion about $\Ad^{r}_1$
follows similarly.
\end{proof}

\begin{remark}
\label{opposite groupoid}
The opposite algebra $H^{op}$ is also a quantum groupoid with the 
same coalgebra structure and the antipode $S^{-1}$. Indeed,
\begin{eqnarray*}
S^{-1}(h\2)h\1
&=& S^{-1}(\eps_s(h)) = S^{-1}(1\1) \eps(h1\2) \\
&=& S^{-1}(1\1) \eps(h S^{-1}(1\2)) =\eps(h1\1)1\2 =\eps_t^{op}(h), \\
h\2 S^{-1}(h\1)
&=& S^{-1}(\eps_t(h)) = \eps(1\1h) S^{-1}(1\2) \\
&=& \eps( S^{-1}(1\1)h)S^{-1}(1\2) = 1\1\eps(1\2h) =\eps_s^{op}(h),\\
S^{-1}(h\3) h\2 S^{-1}(h\1) 
&=& S^{-1}(h\1 S(h\2) h\3) = S^{-1}(h). 
\end{eqnarray*}
Similarly, the co-opposite coalgebra
$H^{cop}$ (with the same algebra structure 
as $H$ and the opposite coalgebra structure, and the antipode $S^{-1}$)
and $H^{op/cop}$ (with both opposite algebra and coalgebra structures, 
and the antipode $S$) are quantum groupoids.
\end{remark}
%%%%%%%%%%%%%%%%%%%%%%%%%%%%%%%%%%%%%%%%%%%%%%%%%%%%%%%%%%%%%%%%%%%%%%%%%%%%%%
%%%%%%%%%%%%%%%% Examples of quantum groupoids %%%%%%%%%%%%%%%%%%%%%%%%%%%%%%%
%%%%%%%%%%%%%%%%%%%%%%%%%%%%%%%%%%%%%%%%%%%%%%%%%%%%%%%%%%%%%%%%%%%%%%%%%%%%%%
\section{Examples of quantum groupoids}
\label{examples}
 
\textbf{Groupoid algebras and their duals ([NV1], 2.1.4).}
As group algebras and their duals give the simplest examples of Hopf
algebras, groupoid algebras and their duals provide 
simple examples of quantum groupoids.

Let $G$ be a finite {\em groupoid} (a category with finitely many
morphisms, such that each morphism is invertible). Then the groupoid
algebra $kG$ (generated by morphisms $g\in G$ with the product of 
two morphisms being equal to  their composition
if the latter is defined and $0$ otherwise) 
is a quantum groupoid via :
\begin{equation}
\Delta(g) = g\otimes g,\quad \eps(g) =1,\quad S(g)=g^{-1},\quad g\in G.
\label{groupoid algebra}
\end{equation}
The counital subalgebras of $kG$ are equal to each other and
coincide with the abelian algebra spanned by the identity morphisms :
$(kG)_t = (kG)_s = \span\{gg^{-1}\mid g\in G\}$. The target and source
counital maps are induced by the operations of taking the target (resp.\
source) object of a morphism : 
\begin{equation*}
\eps_t(g) =gg^{-1} = \id_{target(g)} \quad \mbox{ and } \quad
\eps_s(g) = g^{-1}g = \id_{source(g)}.
\end{equation*} 

The dual quantum groupoid $\widehat{kG}$ is
isomorphic to the algebra of functions on $G$, i.e.,
it is generated by idempotents  $p_g,\, g\in G$ such that
$p_g p_h= \delta_{g,h}p_g$, with

\begin{equation}
\Delta(p_g) =\sum_{uv=g}\,p_u\otimes p_v,\quad \eps(p_g)= \delta_{g,gg^{-1}},
\quad S(p_g) =p_{g^{-1}}.
\label{dual groupoid algebra}
\end{equation}
The target (resp.\ source) counital subalgebra is precisely the algebra
of functions constant on each set of morphisms of $G$ having
the same target (resp.\ source). The target and source maps are
\begin{equation*}
\eps_t(p_g) = \sum_{vv^{-1}=g}\, p_v \quad \mbox{ and } \quad 
\eps_s(p_g) = \sum_{v^{-1}v=g}\, p_v.
\end{equation*}

\begin{definition}
\label{semisimple}
We call a quantum groupoid {\em semisimple} if its underlying algebra
is semisimple.
\end{definition}

In contrast to finite
dimensional semisimple Hopf algebras, the antipode in a finite
dimensional quantum groupoid is not necessarily involutive, see Section 9.

Groupoid algebras and their duals give examples of commutative
and cocommutative semisimple quantum groupoids.
\medskip

%%%%%%%%%%%%%%%%%%%%%%%%%%%%%%%%%%%%%%%%%%%%%%%%%%%%%%%%%%%%%%%%%%%%%%%%%%
\textbf{Quantum transformation groupoids.}
It is known that any group action on a set (i.e., on a commutative algebra 
of functions) gives rise to a groupoid \cite{R}. 
Extending this construction, we associate a quantum groupoid 
with any action of a Hopf algebra on a separable algebra 
(``finite quantum space'').
Namely, let $H$ be a Hopf algebra and $B$ be a separable 
(and, therefore, finite dimensional and  semisimple \cite{P}) 
algebra with right $H$-action $ b\otimes h \mapsto b\cdot h$,
where $b\in B,h\in H$. Then $B^{op}$, the algebra opposite to $B$,
becomes a left $H$-module via $h\otimes a \mapsto 
h\cdot a = a\cdot S_H(h)$.
One can form a {\em double crossed product algebra} $B^{op}\rtimes H\ltimes B$
with underlying vector space $ B^{op}\otimes H \otimes B$ and multiplication
\begin{equation*}
(a\otimes h\otimes b)(a'\otimes h'\otimes b') =
(h\1\cdot a')a \otimes h\2{h'}\1 \otimes (b\cdot {h'}\2) b',
\end{equation*}
for all $a,a'\in B^{op},\, b,b'\in B,$ and $h,h'\in H$.

Assume that $k$ is algebraically closed and let
$e$ be a separability element of $B$ 
(note that $e$  is an idempotent when considered in $B^{op}\otimes B$).
Let $\omega\in \widehat B$ be uniquely determined by 
$(\omega\otimes \id)e = (\id\otimes \omega)e = 1$.
One can check that $\omega$ is the trace of the left
regular representation of $B$ and 
$$
\omega((h\cdot a)b) = \omega(a (b\cdot h)), \qquad
e\I \otimes (h\cdot e\II)  = (e\I \cdot h) \otimes e\II,
$$
where $a\in B^{op}, b\in B$, and $e =e\I\otimes e\II$.

The structure of a quantum groupoid on  $B^{op}\rtimes H\ltimes B$ is given by
\begin{eqnarray} 
\Delta(a\otimes h\otimes b) &=& 
(a\otimes h\1 \otimes e\I)\otimes ((h\2\cdot e\II)\otimes h\3 \otimes b),\\
\eps(a\otimes h\otimes b) &=& \omega(a(h\cdot b)) = \omega(a(b\cdot S_H(h))),\\
S(a\otimes h\otimes b) &=& b\otimes S_H(h)\otimes a.
\end{eqnarray}  
\medskip

%%%%%%%%%%%%%%%%%%%%%%%%%%%%%%%%%%%%%%%%%%%%%%%%%%%%%%%%%%%%%%%%%%%%%%%%%%%
\textbf{Quantum groupoids $\mathbf{B^{op}\otimes B}$ (\cite{BSz2}, 5.2).}
\label{Bop tensor B}
Let $k$ be algebraically closed and let
$B$ be a separable algebra over $k$, $e =e\I \otimes e\II \in
B^{op}\otimes B$ be the symmetric separability idempotent
of $B$, and $\omega$ be as in the previous example.
The map $\pi:x\mapsto e\I x e\II$ defines a linear projection
from $B$ to $Z(B)$. Let $q$ be an invertible element of $B$ such that 
$\pi(q) =1$, then the following operations define a structure
of quantum groupoid $H_q$ on $B^{op}\otimes B$  :  
\begin{eqnarray}
\Delta(b\otimes c) &=& (b\otimes e\I q^{-1}) \otimes (e\II \otimes c), \\
\eps(b\otimes c) &=& \omega(qbc),\\
S(b\otimes c) &=& q^{-1}cq\otimes b,
\end{eqnarray}
for all $b,c\in B$.
The target and source counital subalgebras
of $H_q$ are $B^{op} \otimes 1$ and $1\otimes B$.
The square of the antipode is a conjugation
by $g_q =q\otimes q$. Since $H_q$ with different $q$ are non-isomorphic, 
this example shows that there
can be uncountably many non-isomorphic semisimple quantum
groupoids with the same underlying algebra (for noncommutative $B$).

This example can be also explained in terms of {\em twisting} of quantum
groupoids (see \cite{EN}, \cite{NV4}).
\medskip

%%%%%%%%%%%%%%%%%%%%%%%%%%%%%%%%%%%%%%%%%%%%%%%%%%%%%%%%%%%%%%%%%%%%%%%%%%%
\textbf{Quantum groupoids from subfactors.}
The initial motivation for
studying quantum groupoids in
\cite{NV1}, \cite{NV2}, \cite{N} was their connection
with depth 2 von Neumann subfactors. This connection was first mentioned
in \cite{O} and was also considered in \cite{BNSz}, \cite{BSz1},
\cite{BSz2}, \cite{NSzW}. It was shown in \cite{NV2} that
quantum groupoids naturally arise as non-commutative symmetries
of subfactors, namely if $N\subset M\subset M_1\subset M_2\subset\dots$
is the Jones tower constructed from a finite index, depth $2$
inclusion $N\subset M$ of II$_1$ factors, then $H=M'\cap M_2$
has a canonical structure of a quantum groupoid
acting outerly on $M_1$ such that $M =M_1^H$ and
$M_2 = M_1\rtimes H$. Furthermore $\widehat H=N'\cap M_1$ is a
quantum groupoid dual to $H$.

In \cite{NV3} this result was extended to arbitrary finite depth, via 
a Galois correspondence and it was shown in (\cite{NV3}, 4) that any inclusion
of type II${}_1$ von Neumann factors with finite index and depth
(\cite{GHJ}, 4.1) gives rise to a quantum groupoid and its coideal subalgebra.

We refer the reader to the survey \cite{NV4} (Sections 8,9) and to the Appendix
of \cite{NV3} for the explanation of how quantum groupoids can be constructed
from subfactors.
\medskip

%%%%%%%%%%%%%%%%%%%%%%%%%%%%%%%%%%%%%%%%%%%%%%%%%%%%%%%%%%%%%%%%%%%%%%%%%%%
\textbf{Temperley-Lieb algebras.}
We describe quantum groupoids arising from type $A_n$ subfactors,
whose underlying algebras are {\em Temperley-Lieb algebras} (\cite{GHJ}, 2.1).

Let $k=\mathbb {C}$,  $\lambda^{-1} = 4\cos^2\frac{\pi}{n+3}\ (n\geq 2)$,
and $e_1, e_2,\dots$ be a sequence of idempotents satisfying,
for all $i$ and $j$, the relations
\begin{eqnarray*}
e_i e_{i\pm 1} e_i &=& \lambda e_i, \\
e_i e_j &=& e_j e_i, \quad \mbox{if } |i-j| \geq 2.
\end{eqnarray*}
Let $A_{k,l}$ be the algebra generated by
$1, e_k, e_{k+1},\dots e_l$ ($k\leq l$), $\sigma$ be the algebra 
anti-automorphism of $H= A_{1,2n-1}$ determined by $\sigma(e_i) = e_{2n-i}$
and $P_{k}\in A_{2n-k, 2n-1} \otimes A_{1,k}$ be the image of the 
separability idempotent of $A_{1,k}$ under $\sigma\otimes \id$.

We denote by $\tau$ the non-degenerate Markov trace (\cite{GHJ}, 2.1)
on $H$ and by $w$ the index of the restriction of $\tau$ on
$A_{n+1, 2n-1}\subset H$, i.e., the unique central element in $A_{n+1, 2n-1}$
such that $\tau(w\,\cdot )$ is equal to the trace of the left regular 
representation of  $A_{n+1, 2n-1}$ (see  \cite{W}).

Then the  following operations define a quantum groupoid structure on $H$ :
\begin{eqnarray*}
\Delta(yz) &=& (z\otimes y)P_{n-1}, \qquad y\in A_{n+1,2n-1},\quad
               z\in A_{1,n-1}\\
\Delta(e_n) &=& (1\otimes w) P_{n} (1\otimes w^{-1}), \\
S(h)      &=& w^{-1}\sigma(h)w, \\
\eps(h)   &=& \lambda^{-n}\tau(hfw), \quad h\in A,  
\end{eqnarray*}
where in the last line
$$
f= \lambda^{n(n-1)/2}(e_ne_{n-1}\cdots e_1)(e_{n+1}e_n\cdots e_2)\cdots
   (e_{2n-1}e_{2n-2}\cdots e_n)
$$
is the Jones projection corresponding to the $n$-step basic construction.

The source and target counital subalgebras of $H= A_{1,2n-1}$  are 
$H_s=A_{n+1, 2n-1}$ and $H_t=A_{1,n-1}$. 
The example corresponding to $n=2$ is a quantum 
groupoid of dimension $13$ with antipode of infinite order
(cf. \cite{NV2}, 7.3).

\end{section}

%%%%%%%%%%%%%%%%%%%%%%%%%%%%%%%%%%%%%%%%%%%%%%%%%%%%%%%%%%%%%%%%%%%%%%
%%%%%%  Representation category of a quantum groupoid  %%%%%%%%%%%%%%%
%%%%%%%%%%%%%%%%%%%%%%%%%%%%%%%%%%%%%%%%%%%%%%%%%%%%%%%%%%%%%%%%%%%%%%

\begin{section}
{Representation category of a quantum groupoid}
Throughout this paper we refer to \cite{T2} for definitions 
related to categories.

For a quantum groupoid $H$
let $\Rep(H)$ be the {\em category of representations} of $H$, 
whose objects are finite dimensional left $H$-modules and whose morphisms are 
$H$-linear homomorphisms. We shall show that 
$\Rep(H)$ has a natural  structure of a monoidal category with duality.

For objects $V,W$ of $\Rep(H)$ set
\begin{equation}
V\otimes W = \{x \in V\otimes_k W \mid x= \Delta(1)\cdot x\}\subset
V\otimes_k W\},
\end{equation}
with the obvious action of $H$ via the comultiplication $\Delta$ (here
$\otimes_k$ denotes the usual tensor product of vector spaces).
Note that $\Delta(1)$ is an idempotent and therefore 
$V\otimes W=\Delta(1)\cdot (V\otimes_k W)$.
The tensor product of morphisms is the restriction of usual tensor 
product of homomorphisms. The standard associativity isomorphisms
$(U\otimes V)\otimes W \to U\otimes (V\otimes W)$
are functorial and satisfy the pentagon condition, since $\Delta$ is
coassociative. We will suppress these isomorphisms and write simply
$U\otimes V\otimes W$.

The target counital subalgebra $H_t\subset H$ has an $H$-module structure  
given by $h\cdot z = \eps_t(hz)$, where $h\in H,\, z\in H_t$.

\begin{lemma}
$H_t$ is the unit object of $\Rep(H)$.
\end{lemma}
\begin{proof}
Define a $k$-linear homomorphism $l_V : H_t\otimes V \to V$ by
$$
l_V( 1\1\cdot z \otimes 1\2\cdot v) = z\cdot v, \qquad z\in H_t,\, v\in V.
$$
This map is $H$-linear, since
\begin{eqnarray*}
l_V( h\cdot(1\1\cdot z \otimes 1\2\cdot v))
&=& l_V( h\1\cdot z\otimes h\2\cdot v) \\
&=& \eps_t(h\1z)h\2 \cdot v = hz\cdot v \\
&=& h\cdot l_V( 1\1\cdot z\otimes 1\2\cdot v),
\end{eqnarray*}
for all $h\in H$. The inverse map $l_V^{-1} : V \to H_t\otimes V$ is given by
$$
l_V^{-1}(v) = S(1\1) \otimes (1\2 \cdot v) = (1\1\cdot 1) \otimes (1\2\cdot v).
$$
The collection
$\{l_V\}_V$ gives a natural equivalence between the functor $H_t\otimes (\ )$
and the identity functor. 
Indeed, for any $H$-linear homomorphism
$f: V\to U$ we have :
\begin{eqnarray*}
l_U \circ (\id\otimes f)(1\1\cdot z\otimes 1\2\cdot v)
&=& l_U(1\1\cdot z\otimes 1\2\cdot f(v)) \\
&=& z\cdot f(v) =f(z\cdot v) \\
&=& f\circ l_V(1\1\cdot z\otimes 1\2\cdot v)
\end{eqnarray*}
Similarly, the $k$-linear homomorphism
$r_V : V\otimes H_t \to V$ defined by
$$
r_V( 1\1\cdot v \otimes 1\2\cdot z) = S(z)\cdot v, \qquad z\in H_t,\, v\in V,
$$
has the inverse $r_V^{-1}(v) = 1\1 \cdot v \otimes 1\2$
and satisfies the necessary properties.

Finally, we can check the triangle axiom 
$id_V\otimes l_W = r_V\otimes id_W:V\otimes H_t\otimes W\to V\otimes W$ for all 
objects $V,W$
of $\Rep(H)$. For $v\in V,\,w\in W$  we have
\begin{eqnarray*}
\lefteqn{
(\id_V\otimes l_W)(1\1\cdot v \otimes 1'\1 1\2\cdot z \otimes 1'\2\cdot w)=}\\
&=& 1\1\cdot v \otimes 1\2z \cdot w\\
&=& 1\1S(z)\cdot v\otimes 1\2\cdot w \\
&=& (r_V\otimes \id_W)(1'\1\cdot v \otimes 1'\2 1\1 \cdot z\otimes 1\2\cdot w),  
\end{eqnarray*}  
therefore, $\id_V\otimes l_W = r_V\otimes \id_W$.
\end{proof}

Using the antipode $S$ of $H$, we can provide $\Rep(H)$ with a duality.
For any object $V$ of $\Rep(H)$, define the action of $H$ on $V^*=
\Hom_k(V,\,k)$ by 
\begin{equation}
\label{dual object}
(h\cdot \phi)(v) = \phi(S(h)\cdot v), 
\end{equation}
where $h\in H, v\in V, \phi\in V^*$. For any morphism 
$f: V\to W$, let $f^* : W^*\to V^*$ be the morphism dual to $f$
(see \cite{T2}, I.1.8).

For any $V$ in $\Rep(H)$, we define the duality morphisms
\begin{equation*}
d_V : V^*\otimes V \to H_t, \qquad b_V: H_t \to V \otimes V^*
\end{equation*}
as follows. For $\sum_j\, \phi^j \otimes v_j\in V^* \otimes V$, set
\begin{equation}
\label{dV}
d_V(\sum_j\, \phi^j \otimes v_j)= \sum_j\, \phi^j(1\1\cdot v_j)1\2.
\end{equation}
Let $\{f_i\}_i$ and $\{\xi^i\}_i$ be bases of $V$ and $V^*$ respectively,
dual to each other. The element $\sum_i\,f_i\otimes \xi^i$ does not
depend on choice of these bases; moreover, for all $v\in V, \phi\in V^*$
one has $\phi = \sum_i\,\phi(f_i)\xi^i$ and $v =\sum_i\,f_i\xi^i(v)$.
Set
\begin{equation}
\label{bV}
b_V(z) =  z\cdot (\sum_i\, f_i \otimes \xi^i).
\end{equation}

\begin{proposition}
\label{monoidal with duality}
The category $\Rep(H)$ is a monoidal category with duality.
\end{proposition}
\begin{proof}
We know already that $\Rep(H)$ is monoidal,  it remains to prove
that $d_V$ and $b_V$ are $H$-linear and satisfy the identities
$$
(\id_V\otimes d_V)(b_V \otimes \id_V) = \id_V, \qquad
(d_V \otimes \id_{V^*})(\id_{V^*}\otimes b_V) = \id_{V^*}.
$$ 
Take $\sum_j\, \phi^j \otimes v_j\in V^*\otimes V, z\in H_t,h\in H$.
Using the axioms of a quantum groupoid, we have
\begin{eqnarray*}
h\cdot d_V(\sum_j\, \phi^j \otimes v_j)
&=& \sum_j\,\phi^j(1\1\cdot v) \eps_t(h1\2) \\
&=& \sum_j\,\phi^j(\eps_s(1\1h)\cdot v_j) 1\2 \\
&=& \sum_j\,\phi^j(S(h\1)1\1h\2 \cdot v_j) 1\2 \\
&=& \sum_j\,(h\1\cdot \phi^j)(1\1\cdot(h\2\cdot v_j))1\2 \\
&=& \sum_j\, d_V(h\1\cdot \phi^j \otimes h\2\cdot v_j) \\
&=& d_V (h\cdot \sum_j\, \phi^j \otimes v_j),
\end{eqnarray*}
therefore, $d_V$ is $H$-linear. To check the $H$-linearity of 
$b_V$ we have to show that $h\cdot b_V(z) = b_V(h\cdot z)$, i.e.,\
that
$$
\sum_i\,h\1z\cdot f_i\otimes h\2 \cdot \xi^i
= \sum_i\,1\1\eps_t(hz)\cdot f_i \otimes 1\2 \cdot \xi^i.
$$
Since both sides of the above equality are elements of $V\otimes_k V^*$,
evaluating the second factor on $v\in V$, we get the equivalent condition
$$
h\1zS(h\2)\cdot v = 1\1\eps_t(hz)S(1\2)\cdot v,   
$$
which is easy to check. Thus, $b_V$ is $H$-linear.

Using the isomorphisms $l_V$ and $r_V$ identifying $H_t\otimes V$,
$V\otimes H_t$, and $V$, for all $v\in V$ and $\phi\in V^*$ we have:
\begin{eqnarray*}
(\id_V\otimes d_V)(b_V \otimes \id_V)(v)
&=& (\id_V\otimes d_V)(b_V(1\1\cdot 1)\otimes 1\2\cdot v) \\
&=& (\id_V\otimes d_V)(b_V(1\2)\otimes S^{-1}(1\1)\cdot v) \\
&=& \sum_i\,(\id_V\otimes d_V)
            (1\2\cdot f_i \otimes 1\3\cdot \xi^i \otimes S^{-1}(1\1)\cdot v) \\
&=& \sum_i\,1\2\cdot f_i \otimes (1\3\cdot \xi^i)(1\1'S^{-1}(1\1)\cdot v)1\2'\\
&=& 1\2 S(1\3) 1\1' S^{-1}(1\1)\cdot v \otimes 1\2' =v, \\
(d_V \otimes \id_{V^*})(\id_{V^*}\otimes b_V)(\phi)
&=& (d_V \otimes \id_{V^*})(1\1\cdot \phi \otimes b_V(1\2)) \\
&=& \sum_i\, (d_V \otimes \id_{V^*})
             (1\1\cdot\phi \otimes 1\2\cdot f_i \otimes 1\3\cdot \xi^i) \\
&=& \sum_i\,(1\1\cdot\phi)(1\1'1\2\cdot f_i)1\2' \otimes 1\3\cdot \xi^i\\
&=& 1\2' \otimes 1\3 1\1 S(1\1'1\2)\cdot \phi = \phi,
\end{eqnarray*}
which completes the proof.
\end{proof}
\begin{remark}
Similarly to the construction of $\Rep(H)$, one can construct a category of
right $H$-modules, in which $H_s$ plays the role of the unit object.
\end{remark}
\end{section}

%%%%%%%%%%%%%%%%%%%%%%%%%%%%%%%%%%%%%%%%%%%%%%%%%%%%%%%%%%%%%%%%%%%%%%%%%%%%
%%%  Quasitriangular quantum groupoids %%%%%%%%%%%%%%%%%%%%%%%%%%%%%%%%%%%%%
%%%%%%%%%%%%%%%%%%%%%%%%%%%%%%%%%%%%%%%%%%%%%%%%%%%%%%%%%%%%%%%%%%%%%%%%%%%%

\begin{section}{Quasitriangular quantum groupoids}
\begin{definition}
\label{QT WHA}
A quasitriangular quantum groupoid is a pair 
($H,\, \R$) where $H$ is a quantum groupoid and
$\R\in \Delta^{op}(1)(H\otimes_k H)\Delta(1)$
satisfying the following conditions :
\begin{equation}
\Delta^{op}(h)\R = \R\Delta(h),
\end{equation}
for all $h\in H$,
where $\Delta^{op}$ denotes the comultiplication opposite to $\Delta$,
\begin{equation}
(\id \otimes \Delta)\R = \R_{13}\R_{12}, \qquad
(\Delta \otimes \id)\R = \R_{13}\R_{23},
\end{equation}
where $\R_{12} = \R\otimes 1$, $\R_{23} = 1\otimes \R$, etc.\ as usual,
and such that there exists $\bR\in \Delta(1)(H\otimes_k H)\Delta^{op}(1)$
with
\begin{equation}
\R\bR = \Delta^{op}(1), \qquad \bR\R = \Delta(1).
\end{equation}
\end{definition}

Note that $\bR$ is uniquely determined by $\R$: if $\bR$ and $\bR'$ are two
elements of $\Delta(1)(H\otimes_k H)\Delta^{op}(1)$ satisfying the previous
equation, then 
$$
\bR=\bR\Delta^{op}(1)=\bR\R\bR'=\Delta(1)\bR'=\bR'.
$$
For any two objects $V$ and $W$ of $\Rep(H)$ define
$c_{V,W}: V\otimes W \to W\otimes V$ as the action of $\R_{21}$ :
\begin{eqnarray}
c_{V,W}(x) = \R\II\cdot x\II \otimes \R\I\cdot x\I,
\end{eqnarray}
where  $x =x\I\otimes x\II\in V\otimes W$, $\R=\R\I\otimes \R\II\in
\Delta^{op}(1)(H\otimes_k H)\Delta(1)$.

\begin{proposition}
\label{braiding}
The family of homomorphisms $\{c_{V,W}\}_{V,W}$ defines a braiding in
$\Rep(H)$. Conversely, if $\Rep(H)$ is braided, then there exists
$\R\in \Delta^{op}(1)(H\otimes_k H)\Delta(1)$,
satisfying the properties of Definition~\ref{QT WHA} and inducing
the given braiding.
\end{proposition}
\begin{proof}   
Note that $c_{V,W}$ is well-defined, since $\R_{21}= \Delta(1)\R_{21}$.
To prove the $H$-linearity of $c_{V,W}$ we observe
that
\begin{eqnarray*}
c_{V,W}(h\cdot x)
&=& \R\II h\2\cdot x\II \otimes \R\I h\1\cdot x\I \\
&=& h\1 \R\II \cdot x\II \otimes h\2 \R\I \cdot x\I =h\cdot (c_{V,W}(x)).
\end{eqnarray*}
The inverse of $c_{V,W}$ is given by
$$
c_{V,W}^{-1}(y) = \bR^{(1)}\cdot y\II \otimes \bR^{(2)}\cdot y\I,
$$
where $y =y\I\otimes y\II\in W\otimes V, \bR=\bR\I\otimes \bR\II$.
Therefore, $c_{V,W}$ is an isomorphism.
Finally, one can verify that the braiding identities
$$
(\id_V \otimes c_{U,W})(c_{U,V}\otimes id_W)= c_{U, V\otimes W}, \qquad
(c_{U,W}\otimes\id_V)(\id_U\otimes c_{V,W}) = c_{U\otimes V,W}
$$
are equivalent to the relations of Definition~\ref{QT WHA}, exactly in
the same way as in the case of Hopf algebras 
(see, for instance, (\cite{T2}, XI, 2.3.1)).
\end{proof}

\begin{lemma}
\label{yz} 
Let $(H,\R)$ be a quasitriangular quantum groupoid.
Then for all $y\in H_s, z\in H_t$ the following six identities hold:
\begin{eqnarray*}
(1\otimes z)\R &=& \R(z\otimes 1), \qquad (y\otimes 1)\R = \R(1\otimes y),\\
(z\otimes 1)\R &=& (1\otimes S(z))\R, \qquad (1\otimes y)\R = (S(y)\otimes 1)\R,\\
\R(y\otimes 1) &=& \R(1\otimes S(y)), \qquad \R(1\otimes z) = \R(S(z)\otimes 1).
\end{eqnarray*}
\end{lemma}
\begin{proof}
Since we have $\Delta^{op}(1)\R =\R = \R\Delta(1)$, the first line is
a consequence of the relation $\Delta(yz) = (z\otimes y)\Delta(1)$,
the second line follows from
$(1\otimes z)\Delta(1)= (S(z)\otimes 1)\Delta(1)$ and
$(y\otimes 1)\Delta(1)= (1\otimes S(y))\Delta(1)$. The last two
identities are proven similarly.
\end{proof}
 
\begin{proposition}
\label{QYBE}
Let $(H,\R)$ be a quasitriangular quantum groupoid.
Then $\R$ satisfies the quantum Yang-Baxter equation :
$$
\R_{12}\R_{13}\R_{23} = \R_{23}\R_{13}\R_{12}.
$$
\end{proposition}  
\begin{proof}
It follows from the first two relations of Definition~\ref{QT WHA}, that
\begin{eqnarray*}
\R_{12}\R_{13}\R_{23}=(\id\otimes \Delta^{op})(\R) \R_{23}=
\R_{23} (\id\otimes \Delta)(\R)=\R_{23}\R_{13}\R_{12}.
\end{eqnarray*} 
\end{proof}

\begin{remark}
Let us define two $k$-linear maps $\R_1, \R_2 : \widehat H \to H$ by
$$
\R_1(\phi) =(\id \otimes \phi)(\R), \quad
\R_2(\phi) =(\phi\otimes \id)(\R), \quad\mbox{for}\ \phi\in \widehat H.
$$
Then the condition $(\id \otimes \Delta)\R = \R_{13}\R_{12}$
is equivalent to $\R_1$ being a coalgebra homomorphism and algebra 
anti-homomorphism and the condition $(\Delta \otimes \id)\R = \R_{13}\R_{23}$ 
is equivalent to $\R_2$ being an algebra homomorphism and coalgebra 
anti-homomorphism.  In other words, $\R_1$ and $\R_2$ are
homomorphisms of quantum groupoids $\widehat H \to H^{\op}$ and
$\widehat H \to H^{\cop}$, respectively.
\end{remark}

\begin{proposition}
\label{properties of R}
For any quasitriangular quantum groupoid $(H,\R)$,
we have: 
\begin{eqnarray*}
(\eps_s\otimes \id)(\R) = \Delta(1),&\ &
(\id\otimes \eps_s)(\R) = (S\otimes \id)\Delta^{op}(1), \\
(\eps_t\otimes \id)(\R) = \Delta^{op}(1),&\ &
(\id\otimes \eps_t)(\R) = (S\otimes \id)\Delta(1),
\end{eqnarray*}
$$
(S\otimes \id)(\R) = (\id\otimes S^{-1})(\R) = \bR,\ (S\otimes S)(\R)= \R. 
$$
\end{proposition}
\begin{proof}
First, using the same argument as in \cite{T2}, XI, 2.1.1, we can show that
$(\eps\otimes \id)(\R) = (\id\otimes \eps)(\R) = 1$. Next, using
Lemma~\ref{yz}, we obtain
\begin{eqnarray*}
(\eps_s\otimes \id)(\R) 
&=& 1\1\eps(\R\I 1\2) \otimes \R\II \\
&=& 1\1\eps(\R\I) \otimes 1\2 \R\II = \Delta(1), \\
(\id \otimes \eps_s)(\R)
&=& \R\I \otimes 1\1 \eps(\R\II 1\2) \\
&=& S(1\2)\R\I \otimes 1\1 \eps(\R\II) = (S\otimes \id)\Delta^{op}(1),\\
(\eps_t\otimes \id)(\R)
&=& \eps(1\1 \R\I)1\2 \otimes \R\II  \\
&=& \eps(\R\I) 1\2 \otimes \R\II 1\1 = \Delta^{op}(1), \\
(\id \otimes \eps_t)(\R)
&=& \R\I \otimes \eps(1\1 \R\II) 1\2 \\
&=& S(1\1 )\R\I \otimes \eps(\R\II)1\2 = (S\otimes \id)\Delta(1).
\end{eqnarray*}

Let $m$ denote multiplication $H\otimes_k H\to H$ in $H$. Set
$m_{12}= m\otimes\id: H^{\otimes 3}\to H^{\otimes 2}$ and
$m_{23}= \id\otimes m: H^{\otimes 3}\to H^{\otimes 2}$.
It follows from the above relations that
\begin{eqnarray*}
m_{12}((S\otimes \id \otimes \id)(\Delta\otimes \id)(\R)) 
&=& (\eps_s\otimes \id)(\R) = \Delta(1), \\
m_{12}((\id \otimes S \otimes \id)(\Delta\otimes \id)(\R))    
&=& (\eps_t\otimes \id)(\R) = \Delta^{op}(1), \\
m_{23}((\id \otimes \id \otimes S^{-1})(\Delta^{op}\otimes \id)(\R))   
&=& (\id \otimes S^{-1}\eps_t)(\R) = \Delta^{op}(1), \\
m_{23}((\id \otimes S^{-1} \otimes \id)(\Delta^{op}\otimes \id)(\R))
&=& (\id\otimes S^{-1}\eps_s)(\R) = \Delta(1).
\end{eqnarray*}
On the other hand, Definition~\ref{QT WHA} implies that
\begin{eqnarray*}
\lefteqn{m_{12}((S\otimes \id \otimes \id)(\Delta\otimes \id)(\R)) =}\\
&=& m_{12}((S\otimes \id \otimes \id)(\R_{13}\R_{23}) ) 
    = (S\otimes\id)(\R)\R,\\
\lefteqn{m_{12}((\id \otimes S \otimes \id)(\Delta\otimes \id)(\R))=}\\
&=& m_{12}((\id \otimes S \otimes \id)(\R_{13}\R_{23})) 
    = \R(S\otimes\id)(\R),\\
\lefteqn{m_{23}((\id \otimes \id \otimes S^{-1})(\Delta^{op}\otimes \id)(\R))=}\\
&=& m_{23}((\id \otimes \id \otimes S^{-1})(\R_{12}\R_{13}))
    = \R(\id \otimes S^{-1})(\R),\\
\lefteqn{m_{23}((\id \otimes S^{-1} \otimes \id)(\Delta^{op}\otimes \id)(\R))=}\\
&=& m_{23}((\id \otimes S^{-1} \otimes \id)(\R_{12}\R_{13}))
    = (\id \otimes S^{-1})(\R)\R.
\end{eqnarray*}
Therefore, $(S\otimes \id)(\R) = (\id\otimes S^{-1})(\R) = \bR$
and  $(S\otimes S)(\R)= \R$.
\end{proof}

\begin{proposition}
\label{elements u and v}  
Let $(H,\R)$ be a quasitriangular quantum groupoid.
Then $S^2(h) = uhu^{-1}$ for all $h\in H$, where $u=S(\R\II)\R\I$ is
an invertible element of $H$ such that
$$
u^{-1} = \R\II S^2(\R\I), \quad \Delta(u) = \bR \bR_{21} (u\otimes u).
$$
Likewise, $S^{-2}(h) = vhv^{-1}$, where $v = S(u) = \R\I S(\R\II)$, and
$$
v^{-1} = S^2(\R\I)\R\II, \quad \Delta(v) = \bR \bR_{21} (v\otimes v).
$$
\end{proposition}
\begin{proof}
Note that $S(\R\II)y\R\I = S(y)u$ for all $y\in H_s$, by Lemma~\ref{yz}. 
Hence, we have
\begin{eqnarray*}
S(h\2)uh\1
&=& S(h\2) S(\R\II)\R\I h\1 = S(\R\II h\2)\R\I h\1 \\
&=& S(h\1 \R\II) h\2 \R\I = S(\R\II) \eps_s(h) \R\I \\
&=& S(\eps_s(h))u,
\end{eqnarray*}
for all $h\in H$. Therefore, using the axioms of a quantum groupoid,
we get
\begin{eqnarray*}
uh
&=& S(1\2)u1\1h = S(\eps_t(h\2)u h\1 \\
&=& S(h\2 S(h\3))u h\1 = S^2(h\3) S(h\2)u h\1 \\
&=&  S^2(h\2) S(\eps_s(h\1)) u = S(\eps_s(h\1)S(h\2))u = S^2(h)u.
\end{eqnarray*}
The remaining part of the proof follows the lines of
(\cite{M}, 2.1.8). The results for $v$ can be obtained
by applying the results for $u$ to the quasitriangular quantum groupoid
$(H^{op/cop}, \R)$.
\end{proof}

\begin{definition}
\label{Drinfeld element}
The element $u$ defined in Proposition~\ref{elements u and v} is
called {\em the Drinfeld element} of $H$.
\end{definition} 
  
\begin{corollary}
\label{uv is central}
The element $uv=vu$ is central and
$$
\Delta(uv) = (\bR \bR_{21})^2(uv\otimes uv).
$$
The element $uv^{-1} =vu^{-1}$ is group-like and
$S^4(h)=uv^{-1}hvu^{-1}$ for all $h\in H$.
\end{corollary}   
\begin{proof}
See (\cite{M}, 2.1.9).
\end{proof}

\begin{proposition}
\label{range of F}
Given a quasitriangular quantum groupoid $(H, \R)$, consider a linear map
$F:\widehat H\to H$ given by
\begin{equation}
F(\phi)=(\phi\otimes \id)(\R_{21}\R),\qquad \phi\in \widehat H.
\end{equation}
Then the image of $F$ lies in $C_H(H_s)$, the centralizer of $H_s$.
\end{proposition}
\begin{proof}
Take $y\in H_s$. Then we have
$$
\phi(\R\II {\R}\I) \R\I {\R}\II y =
\phi(\R\II y {\R}\I) \R\I {\R}\II = \phi(\R\II {\R}\I) y \R\I {\R}\II.
$$
Therefore $F(\phi)\in C_H(H_s)$, as required.
\end{proof}

\begin{definition}[cf. \cite{M}, 2.1.12]
\label{factorizability}
A quasitriangular quantum groupoid is {\em factorizable} if
the map   $F:\widehat H \to C_H(H_s)$ from Proposition~\ref{range of F}
is surjective.
\end{definition}

The factorizability means that $\R$
is as non-trivial as possible, in contrast to {\em triangular} 
quantum groupoids, for which $ \bR=\R_{21}$ and $F(\widehat H)=H_t$.
\begin{corollary}
\label{restriction of F}
If $H$ is factorizable, then the restriction of $F$ to the subspace
$W_s = \{ \phi\in \widehat H \mid \phi = \phi\circ\Ad_1^r) \}$
is a linear isomorphism onto $C_H(H_s)$. 
\end{corollary}
\begin{proof}
From the observation that  $F(\phi) = F(\phi\circ\Ad_1^r)$ we have that
the restriction of $F$ to $W_s$ is a linear map onto $C_H(H_s)$.
On the other hand, Lemma~\ref{expectations} allows to identify $W_s$ 
with the dual vector space to $C_H(H_t)$, from where 
$\dim W_s = \dim C_H(H_s)$ and the result follows.
\end{proof}
\end{section}

%%%%%%%%%%%%%%%%%%%%%%%%%%%%%%%%%%%%%%%%%%%%%%%%%%%%%%%%%%%%%
%%%%%%%%%%%%%%%   The Drinfeld double construction %%%% %%%%%
%%%%%%%%%%%%%%%%%%%%%%%%%%%%%%%%%%%%%%%%%%%%%%%%%%%%%%%%%%%%%

\begin{section}{The Drinfeld double for quantum groupoids}

Let $H$ be a finite quantum groupoid.
We define the {\em Drinfeld double} $D(H)$ of $H$ as follows.
Consider on the vector space $\widehat H^{op}\otimes_k H$ a
multiplication given by
\begin{equation}
(\phi\otimes h)(\psi\otimes g) =
\psi\2\phi \otimes h\2g \la S(h\1),\, \psi\1\ra \la h\3,\, \psi\3\ra,
\end{equation}
where $\phi, \psi \in {\widehat H}^{\op}$ and $h,g \in H$. 
We verify below that the linear span $J$ of the elements
\begin{eqnarray}
\label{amalgamation}
\phi\otimes zh &-& (\eps\actl z)\phi \otimes h, \quad z\in H_t,\\
\phi\otimes yh &-& (y\actr\eps)\phi \otimes h, \quad y\in H_s, 
\end{eqnarray}
is a two-sided ideal in $\widehat H^{op}\otimes_k H$. Let $D(H)$ be 
the factor-algebra ($\widehat H^{op}\otimes_k H)/J$ and let 
$[\phi\otimes h]$ denote the class of $\phi\otimes h$ in $D(H)$. 

\begin{definitionandtheorem}
\label{the double} $D(H)$ is a quantum groupoid with unit $[\eps\otimes 1]$,
and comultiplication, counit, and antipode given by 
\begin{eqnarray}
\Delta([\phi\otimes h]) 
&=& [\phi\1 \otimes h\1] \otimes [\phi\2 \otimes h\2],\\
\E([\phi\otimes h])
&=& \la \eps_t(h),\,\phi \ra,\\
S([\phi\otimes h]) 
&=& [ S^{-1}(\phi\2) \otimes S(h\2)] \la h\1,\,\phi\1 \ra
                                     \la S(h\3),\,\phi\3\ra. 
\end{eqnarray}
\end{definitionandtheorem}
In the case where $H$ is a Hopf algebra, this definition and
is due to Drinfeld \cite{D}.
\begin{proof}
Associativity of multiplication in $\widehat H^{op}\otimes_k H$ and hence in
$D(H)$ can be verified exactly as in (\cite{M}, 7.1.1). 
Let us check that $J$ is an ideal. We have :
\begin{eqnarray*}
(\phi\otimes h)((\eps\actl z)\psi \otimes g)
&=& \psi\2\phi \otimes h\2g   
    \la S(h\1),\, (\eps\actl z)\psi\1\ra   \la h\3,\, \psi\3\ra   \\
&=& \psi\2\phi \otimes h\3g   \la zS(h\2),\,\eps\ra 
    \la S(h\1),\, \psi\1\ra   \la h\4,\, \psi\3\ra          \\
&=& \psi\2\phi \otimes h\2zg   \la S(h\1),\, \psi\1\ra  
    \la h\3,\, \psi\3\ra   \\
&=& (\phi\otimes h)(\psi\otimes zg), \\
(\psi\otimes zg)(\phi\otimes h)
&=& \phi\2\psi \otimes g\2h   \la S(zg\1),\,\phi\1 \ra
    \la g\3,\, \phi\3\ra   \\
&=& \la Sz,\, \phi\2\ra \phi\3 \psi \otimes g\2h   \la S(g\1),\,\phi\1 \ra
    \la g\3,\, \phi\3\ra   \\                                
&=& \phi\2 (\eps\actl z)\psi \otimes  g\2h   \la S(g\1),\,\phi\1 \ra 
    \la g\3,\, \phi\3\ra   \\
&=& ( (\eps\actl z)\psi \otimes g)(\phi\otimes h),  
\end{eqnarray*}
where $z\in H_t$ and we used the identity
$zS(h\1)\otimes h\2 = S(h\1)\otimes h\2z$.
Similarly, one checks that
\begin{eqnarray*}
(\psi\otimes yg)(\phi\otimes h) 
&=& ((y\actr \eps)\psi \otimes g)(\phi\otimes h), \\
(\phi\otimes h)(\psi\otimes yg)
&=& (\phi\otimes h) ((y\actr \eps)\psi \otimes g),
\end{eqnarray*}
therefore for all $x\in J$ we have $(\phi\otimes h)x = x(\phi\otimes h) =0$,
so $J$ is an ideal. We also compute
\begin{eqnarray*}
[\eps\otimes 1][\phi\otimes h]
&=& [\phi\2 \otimes 1\2h \la S(1\1),\, \phi\1\ra \la 1\3,\, \phi\3\ra] \\
&=& [\phi\2 \otimes  \la S(1\1),\, \phi\1\ra 1\21'\1 \la 1'\2,\, \phi\3\ra h]\\
&=& [\eps_t(\phi\1)S(\eps_t(\phi\3))\phi\2 \otimes h] = [\phi\otimes h], 
\end{eqnarray*}
and similarly $[\phi\otimes h][\eps\otimes 1] = [\phi\otimes h]$, so that
$[\eps\otimes 1]$ is a unit.

Now let us verify that the structure maps $\Delta,\, \E,$ and $S$
are well-defined on $D(H)$. We have, using properties of a quantum
groupoid and its counital subalgebras:
\begin{eqnarray*}
\Delta([\phi\otimes zh])
&=& [\phi\1\otimes zh\1] \otimes [\phi\2\otimes h\2] \\
&=& [(\eps\actl z)\phi\1 \otimes h\1] \otimes [\phi\2\otimes h\2] \\
&=& \Delta([ \la z,\,\eps\1\ra \eps\2\phi  \otimes h]), \\
\eps([ \la z,\,\eps\1\ra \eps\2\phi  \otimes h])
&=& \la z,\, \eps\1 \ra  \la \eps_t(h),\, \eps\2\phi \ra \\
&=& \la z,\, \eps\1 \ra \la 1\1\eps_t(h),\,\eps\2 \ra \la 1\2,\, \phi\ra \\
&=& \la z\eps_t(h)1\1,\, \eps \ra   \la 1\2,\, \phi\ra \\
&=& \la z\eps_t(h),\,\phi\ra = \E([\phi\otimes zh]), \\
S([ (\eps\actl z)\phi  \otimes h])
&=&  [ \la z,\, \eps\1\ra S^{-1}(\eps\3\phi\2) \otimes S(h\2)] \\ 
& &  \la h\1,\,\eps\2\phi\1 \ra   \la S(h\3),\,\eps\4\phi\3 \ra \\
&=&  [ S^{-1}(\phi\2) \otimes S(h\2)] \\
& &  \la h\1,\,(\eps\actl z)\phi\1 \ra   \la S(h\3),\,\phi\3 \ra \\
&=&  [S^{-1}(\phi\2) \otimes S(h\3) ]\la zh\1,\, \eps\ra \\
& &  \la h\2,\, \phi\1\ra \la S(h\4),\,\phi\3 \ra \\ 
&=&  [S^{-1}(\phi\2) \otimes S(h\2)]  \la zh\1,\, \phi\1 \ra
     \la S(h\3),\,\phi\3 \ra  \\
&=& S([\phi\otimes zh]) 
\end{eqnarray*}
for all $h\in H, \phi\in \widehat{H}, z\in H_t$.
Next, we need to check the axioms of a quantum groupoid. Coassociativity
and multiplicativity of $\Delta$ are established as in (\cite{M}, 7.1.1),
since the computations given there do not use the unitality of multiplication
and comultiplication. For the counit property,
we have:
\begin{eqnarray*}
(\E\otimes\id)\Delta([\phi\otimes h])
&=& \la \eps_t(h\1),\, \phi\1 \ra [\phi\2\otimes h\2] \\
&=& \la \eps_t(h\1),\, \eps\1 \ra [\eps\2\phi \otimes h\2] \\
&=& [\phi \otimes \eps_t(h\1)h\2] = [\phi\otimes h], \\
(\id \otimes\E)\Delta([\phi\otimes h])
&=& [\phi\1\otimes h\1] \la \eps_t(h\2),\, \phi\2 \ra \\
&=& [\phi\1\otimes 1\1h] \la 1\2,\, \eps_t(\phi\2) \ra \\
&=& [S^{-1}(\eps_t(\phi\2))\phi\1 \otimes h] = [\phi\otimes h],
\end{eqnarray*}
where  we used the amalgamation property 
$[\phi\otimes zh]=[(\eps\actl z)\phi \otimes h], \quad z\in H_t,$ 
following from (\ref{amalgamation}).
Now we verify the remaining axioms of a quantum groupoid.
For all $h,g,f\in H$ and $\phi,\psi,\theta\in \widehat{H}$
we compute
\begin{eqnarray*}
\lefteqn{\eps([\phi\otimes h][\psi\otimes g][\theta\otimes f]) = } \\
&=& \eps([\theta\2 \psi\2 \phi \otimes h\3g\2f]) \la S(h\1),\, \psi\1\ra \\
& & \la S(h\2g\1),\theta\1\ra   \la h\4g\3,\, \theta\3 \ra 
    \la h\5,\, \psi\3 \ra \\
&=& \la h\3g\2,\,\eps\1\eps_t(\theta\2 \psi\2 \phi)\ra \la f,\,\eps\2 \ra
    \la S(h\1),\, \psi\1\ra \\
& & \la h\2g\1,\, S(\theta\1)\ra   \la h\4g\3,\, \theta\3 \ra
    \la h\5,\, \psi\3 \ra \\
&=& \la h\2g,\, S(\theta\1)\eps\1 \eps_t(\theta\2 \psi\2 \phi) \theta\3 \ra
    \la f,\,\eps\2 \ra   \la S(h\1),\, \psi\1\ra  \la h\3,\, \psi\3 \ra \\
&=& \la h\2g,\, \eps_t(\psi\2 \phi)\eps_s(\eps\1\theta) \ra
    \la f,\,\eps\2 \ra   \la S(h\1),\, \psi\1\ra  \la h\3,\, \psi\3 \ra \\
&=& \la h\2g\1,\, \eps_t(\psi\2 \phi) \ra  
    \la g\2,\, \eps'\2 \eps_s(\eps\1\theta) \ra \\
& & \la f,\,\eps\2 \ra   \la S(h\1),\, \psi\1\ra  \la h\3,\, \eps'\1 \psi\3\ra\\
&=& \la h\2g\1,\, \eps_t(\psi\2 \phi) \ra 
    \la g\2,\, S(\theta\1)\eps\1 \eps_t(\theta\2 \psi\4)\theta\3 \ra \\
& & \la f,\,\eps\2 \ra  \la S(h\1),\, \psi\1\ra \la h\3,\, \eps'\1 \psi\3 \ra\\
&=& \la h\2g\1,\, \eps_t(\psi\2 \phi) \ra  
     \la S(h\1),\, \psi\1\ra  \la h\3,\, \psi\3 \ra \\
& & \la g\3f,\,\eps_t(\theta\2 \psi\4)\ra 
    \la S(g\2),\, \theta\1 \ra \la g\4,\, \theta\3 \ra \\
&=&  \eps([\phi\otimes h][\psi\1\otimes g\1])
      \eps([\psi\2\otimes g\2][\theta\otimes f]), \\
\lefteqn{\eps([\phi\otimes h][\psi\2\otimes g\2])             
      \eps([\psi\1\otimes g\1][\theta\otimes f]) =} \\
&=& \la h\2 g\4,\, \eps_t(\psi\3 \phi) \ra
     \la S(h\1),\, \psi\2\ra  \la h\3,\, \psi\4 \ra \\
& & \la g\2f,\,\eps_t(\theta\2 \psi\1)\ra
    \la S(g\1),\, \theta\1 \ra \la g\3,\, \theta\3 \ra \\
&=& \la h\2 g\2,\,\eps_t(\psi\3 \phi) \ra
    \la g\1,\, S(\theta\1)\eps\1 \eps_t(\theta\2 \psi\1) \theta\3 \ra \\
& & \la f,\, \eps\2\ra  \la S(h\1),\, \psi\2\ra  \la h\3,\, \psi\4 \ra \\
&=& \la h\2 g\2, \,\eps_t(\psi\3 \phi) \ra
     \la g\1,\, \eps'\2\eps_s(\eps\1\theta) \ra \\
& & \la f,\, \eps\2\ra \la S(h\1),\, S(\eps'\1)\psi\1\ra  \la h\3,\, \psi\4 \ra \\
&=& \la h\2g,\, \eps_t(\psi\2 \phi)\eps_s(\eps\1\theta) \ra
    \la f,\,\eps\2 \ra   \la S(h\1),\, \psi\1\ra  \la h\3,\, \psi\3 \ra \\
&=& \eps([\phi\otimes h][\psi\otimes g][\theta\otimes f]), 
\end{eqnarray*}
which is axiom (\ref{eps m}). For axiom (\ref{Delta 1}) we have:
\begin{eqnarray*}
\lefteqn{(\Delta([\eps\otimes 1])\otimes [\eps\otimes 1])
([\eps\otimes 1] \otimes \Delta([\eps\otimes 1]) ) =}\\
&=& [\eps\1\otimes 1\1] \otimes [\eps\2\otimes 1\2][\eps\1'\otimes 1\1']
    \otimes [\eps\2'\otimes 1\2'] \\
&=& [\eps\1\otimes 1\1] \otimes [\eps\1' \eps\2 \otimes 1\2 1\1']
    \otimes [\eps\2'\otimes 1\2'] \\
&=&  [\eps\1\otimes 1\1] \otimes  [\eps\2\otimes 1\2]
    \otimes  [\eps\3\otimes 1\3], \\
\lefteqn{([\eps\otimes 1] \otimes \Delta([\eps\otimes 1]) )
(\Delta([\eps\otimes 1])\otimes [\eps\otimes 1]) =}\\
&=& [\eps\1'\otimes 1\1']  \otimes  [\eps\1 \otimes 1\1] [\eps\2'\otimes 1\2']
    \otimes  [\eps\2\otimes 1\2] \\
&=& [\eps\1'\otimes 1\1']  \otimes  [\eps\2'\eps\1 \otimes 1\11\2']
   \otimes  [\eps\2\otimes 1\2] \\
&=&  [\eps\1\otimes 1\1] \otimes  [\eps\2\otimes 1\2]
    \otimes  [\eps\3\otimes 1\3], 
\end{eqnarray*}
where we used the axioms of a quantum groupoid and
the definition of $J$. 

In order to check the axioms (\ref{S epst}),(\ref{S epss}), let us compute
the  target counital map $\E_t$. We have
\begin{eqnarray*}
\E_t([\phi\otimes h])
&=& \E([\eps\1\otimes 1\1][\phi\otimes h]) [\eps\2\otimes 1\2] \\
&=& \la \eps_t(1\2h),\,\phi\2\eps\1 \ra   \la S(1\1),\,\phi\1 \ra
    \la 1\3,\, \phi\3 \ra   [\eps\2\otimes 1\4] \\
&=& \la 1'\1\eps_t(1\2h),\,\phi\2\ra  \la 1'\2,\,\eps\1  \ra \\
& & \la S(1\1),\,\phi\1 \ra \la 1\3,\, \phi\3 \ra   [\eps\2\otimes 1\4] \\
&=& \la S(1\1)1'\1 \eps_t(1\2h) 1\3,\, \phi \ra   \la 1'\2,\,\eps\1  \ra
     [\eps\2\otimes 1\4] \\
&=& \la 1\1 \eps_t(h),\, \phi\ra [\eps \otimes 1\2].
\end{eqnarray*}
Similarly one computes the source counital map :
$$
\E_s([\phi\otimes h]) = 
[\eps\1\otimes 1] \la h,\, \eps_t(\phi)S(\eps\2) \ra.
$$
Using these formulas we have :
\begin{eqnarray*}
\lefteqn{m(\id\otimes S)\Delta([\phi\otimes h])= } \\
&=& [\phi\1\otimes h\1][S^{-1}(\phi\3)\otimes S(h\2)]
    \la h\2,\, \phi\2 \ra   \la S(h\4),\, \phi\4 \ra  \\
&=& [ S^{-1}(\phi\4)\phi\1 \otimes h\2 S(h\5)]\la S(h\1),\,S^{-1}(\phi\5) \ra\\
& & \la h\3 ,\,S^{-1}(\phi\3) \ra  \la h\4,\, \phi\2 \ra
    \la S(h\6) ,\, \phi\6 \ra   \\
&=& [ S^{-1}(\phi\3)\phi\1 \otimes h\2 S(h\4)] 
    \la \eps_t(h\3),\, S^{-1}(\phi\2)\ra \la  h\1 S(h\5),\,\phi\4\ra   \\
&=& [ S^{-1}(\phi\3)\phi\1 \otimes 1\1\eps_t(h\2) ]
    \la 1\2,\, S^{-1}(\phi\2)\ra  \la  h\1 S(h\3),\,\phi\4\ra   \\
&=& [ S^{-1}(\phi\3)\phi\1 \otimes 1\1 1\2']
    \la 1\2,\, S^{-1}(\phi\2)\ra  \la  1\1'\eps_t(h),\,\phi\4\ra   \\
&=& [ \eps\1 S^{-1}(\phi\3)\phi\1 \otimes 1\1 1\2'] 
    \la 1\2,\,\eps\2\ra \la  1\1'\eps_t(h),\,\phi\3 \ra   \\
&=& [ S^{-1}(\phi\2)\phi\1 \otimes 1\2']  \la  1\1'\eps_t(h),\,\phi\3 \ra \\
&=& [\eps\otimes 1\2] \la 1\1\eps_t(h),\,\phi\ra = \eps_t([\phi\otimes h]),
\end{eqnarray*}
and
\begin{eqnarray*}
\lefteqn{m(S\otimes\id)\Delta([\phi\otimes h])= } \\
&=& [\phi\5 S^{-1}(\phi\2) \otimes S(h\3)h\6] \la S^{2}(h\4),\,\phi\4\ra
    \la S(h\2),\, \phi\6 \ra \\
& &  \la h\1,\, \phi\1\ra \la S(h\5),\,\phi\3\ra \\
&=& [ S^{-1}(\eps\1\eps_t(\phi\2))\otimes S(h\2) h\4 ]
    \la S(h\3),\,\eps\2\ra  \la h\1,\, \phi\1S(\phi\3) \ra \\
&=& [ S^{-1}(\eps\1\eps\2') \otimes S(h\2)h\4]
    \la \eps_s(h\3),\, S(\eps\2) \ra  \la h\1 \eps\1'\eps_t(\phi) \ra \\
&=& [ S^{-1}(\eps\1\eps\2') \otimes  \eps_s(h\2)]
    \la 1\1,\,  S(\eps\2) \ra  \la h\1 \eps\1'\eps_t(\phi) \ra \\
&=& [\eps\3S^{-1}(\eps\2)  \otimes  \eps_s(h\2)]
    \la h\1,\, \eps_t(\phi)S(\eps\2) \ra \\
&=& [\eps\1 \otimes 1] \la h,\, S(\eps\2)\eps_t(\phi) \ra  =
     \eps_s([\phi\otimes h]).
\end{eqnarray*}
In the above  computations we used repeatedly the amalgamation 
relations in $D(H)$: 
\begin{equation*}
[\phi\otimes zh]=[(\eps\actl z)\phi \otimes h]\ \quad (z\in H_t), \qquad
[\phi\otimes yh]=[(y\actr\eps)\phi \otimes h]\  \quad (y\in H_s)
\end{equation*} 
that follow from (\ref{amalgamation}), the axioms of a quantum groupoid
and properties of  the counital maps. Finally, we prove the relation which is
equivalent to $S$ being both algebra and coalgebra anti-homomorphism:
\begin{eqnarray*}
\lefteqn{m(\id\otimes m)(S\otimes \id\otimes S)
         (\id\otimes\Delta)\Delta([\phi\otimes h])= } \\
&=& m(S\otimes \eps_t)\Delta([\phi\otimes h]) \\
&=& S([\phi\1\otimes h\1])\eps_t([\phi\2\otimes h\2]) \\
&=& [S^{-1}(\phi\2) \otimes S(h\2)] [\eps \otimes 1\2] \\
& &   \la h\1,\, \phi\1 \ra  \la S(h\3),\, \phi\3 \ra
      \la 1\1\eps_t(h\4),\, \phi\4 \ra \\
&=& [S^{-1}(\phi\2) \otimes S(h\2)1\2] \la h\1,\, \phi\1 \ra 
    \la S(h\3)1\1,\,\phi\3 \ra = S([\phi\otimes h]).
\end{eqnarray*}
Note that $D(H)_t = [\eps \otimes H_t]$ and $D(H)_s = [\widehat H_s \otimes 1]$.
\end{proof}

\begin{proposition}
\label{QT of D(H)}
The Drinfeld double $D(H)$ has a quasitriangular structure given by
\begin{equation}
\R = \sum_i [\xi^i \otimes 1] \otimes [\eps \otimes f_i],
\qquad
\bR  = \sum_j\, [S^{-1}(\xi_j) \otimes 1] \otimes[\eps \otimes f_j]
\end{equation}
where $\{ f_i\}$ and $\{ \xi^i\}$ are dual bases in $H$ and $\widehat H$.
\end{proposition}
\begin{proof}
The identities $(\id \otimes \Delta)\R = \R_{13}\R_{12}$ and
$(\Delta \otimes \id)\R = \R_{13}\R_{23}$
can be written as (identifying $[\widehat H^{op}\otimes 1]$ with $\widehat H^{op}$ 
and $[\eps \otimes H]$ with $H$):
\begin{equation}
\label{R-dual bases1}
\sum_i\, \xi\1^i \otimes \xi\2^i \otimes f_i =
\sum_{ij}\, \xi^i \otimes \xi^j \otimes f_i f_j, 
\end{equation}
\begin{equation}
\label{R-dual bases2}
\sum_i\, \xi^i \otimes {f_i}\1  \otimes {f_i}\2 =
\sum_{ij}\, \xi^j \xi^i \otimes f_j \otimes f_i.
\end{equation}
The above equalities can be verified by evaluating both sides
on an element $h\in H$ in the third factor (resp.,\ on
$\phi\in \widehat H^{op}$ in the second factor), see (\cite{M}, 7.1.1).
To show that $\R$ is an intertwiner between $\Delta$ and $\Delta^{op}$,
we compute   
\begin{eqnarray*}
\R \Delta([\phi\otimes h])
&=& \sum_i\, [\phi\1 \xi^i \otimes h\1] \otimes [\phi\3 \otimes
    {f_i}\2 h\2] \\
& & \la S({f_i}\1),\, \phi\2 \ra  \la {f_i}\3,\, \phi\4 \ra \\
&=& \sum_i\, [\phi\1 S(\phi\2)\xi^i\phi\4 \otimes h\1] \otimes
    [\phi\3 \otimes f_i h\2] \\
&=& \sum_i\, [ \xi^i\phi\3 \otimes \la 1\1,\, \eps_t(\phi\1)\ra 1\2h\1 ]
    \otimes [\phi\2 \otimes f_i h\2] \\
&=& \sum_i\, [ \xi^i\phi\2 \otimes \la 1\1,\,\eps\1 \ra 1\2h\1 ]
    \otimes [\eps\2 \phi\1 \otimes f_i h\2] \\
&=& \sum_i\, [ \xi^i\phi\2 \otimes h\2] \otimes
    [ \la \eps_t(h\1),\,\eps\1 \ra \eps\2 \phi\1 \otimes f_i h\3] \\
&=& \sum_i\, [ \xi^i\phi\2 \otimes h\3] \otimes
    [ \phi\1 \otimes h\1 S(h\2) f_i h\4 ] \\
&=& \sum_i\, [\xi\2^i\phi\2 \otimes h\3] \otimes [\phi\1\otimes h\1 f_i]
    \la S(h\2),\, \xi^i\1 \ra   \la h\4,\,\xi^i\3 \ra \\
&=& \Delta^{\op}([\phi\otimes h]) \R.
\end{eqnarray*}
where we used
\begin{equation}
\label{dual bases}
\sum_i\,\la a,\, \xi^i\ra f_i =a\quad \mbox{ and } \quad
\sum_i\, \xi^i \la f_i,\,\phi \ra =\phi,
\end{equation} 
for all $a\in H, \phi\in \widehat H$.
Finally, let us check that the element $\bR = \sum_j\, [S^{-1}(\xi_j)
\otimes 1] \otimes[\eps \otimes f_j]$ satisfies $\bR \R = \Delta(1)$ and
$\R\bR = \Delta^{op}(1)$. The first property is equivalent to
$$
\sum_{i,j}\, [\xi^i S^{-1}(\xi^j)\otimes 1] \otimes 
[\eps \otimes f_jf_i] =
[\la 1\1,\,\eps'\2\ra \eps'\1 \eps\1 \otimes 1] \otimes
[\eps \otimes \la 1'\1,\, \eps\2 \ra 1'\2 1\2],
$$ 
which can be regarded as an equality in $\widehat H^{op}\otimes H$ :
\begin{equation*}
\sum_{i,j}\, \xi^i S^{-1}(\xi^j)\otimes f_jf_i =
\la 1\1,\,\eps'\2\ra \eps'\1 \eps\1 \otimes 
\la 1'\1,\, \eps\2 \ra 1'\2 1\2.
\end{equation*}
Evaluating both sides on arbitrary $\phi\in \widehat H$ in the second factor, 
we get
\begin{eqnarray*}
\phi\2 S^{-1}(\phi\1) 
&=& \la 1\1,\,\eps'\2\ra \eps'\1 \eps\1 \la 1'\1,\, \eps\2 \ra
    \la 1'\2 1\2,\phi \ra \\
&=& \eps_s^{op}(\phi\2) \eps_s^{op}(\phi\1),
\end{eqnarray*}
where $\eps_s^{op}(\phi) = \phi\2 S^{-1}(\phi\1)$ is the source
counital map in ${\widehat H}^{op}$.
The second property is similar.
\end{proof}
For the Hopf algebra case the above idea of the proof was proposed in 
\cite{M}, 7.1.1.

\begin{remark}
The dual quantum groupoid $\widehat{D(H)}$ consists of all elements
$\sum_k\, h_k\otimes \phi_k$ in $H\otimes_k {\widehat H}^{op}$ such that
$$
\sum_k\, (h_k\otimes \phi_k)|_J = 0.
$$
The structure operations in $\widehat{D(H)}$ are obtained by dualizing those
in $D(H)$:
\begin{eqnarray*}
(  \sum_k\, h_k\otimes \phi_k ) (  \sum_l\, g_l\otimes \psi_l )
&=& \sum_{k,l}\, h_kg_l\otimes \phi_k \psi_l, \\
1_{\widehat{D(H)}}
&=& 1\2 \otimes (\eps \actl 1\1), \\
\Delta(\sum_k\, h_k\otimes \phi_k )
&=& \sum_{i,j,k}\, ({h_k}\2 \otimes \xi^i {\phi_k}\1 \xi^j) \otimes
    (S(f_i) {h_k}\1 f_j \otimes {\phi_k}\2), \\
\eps(  \sum_k\, h_k\otimes \phi_k ) 
&=& \sum_k\,\eps(h_k)\widehat\eps(\phi_k),\\
S(  \sum_k\, h_k\otimes \phi_k )
&=& \sum_{i,j,k}\, f_i S^{-1}(h_k)S(f_j) \otimes \xi^i S(\phi)\xi^j,
\end{eqnarray*}
for all $\sum_k\, h_k\otimes \phi_k,\, \sum_l\, g_l\otimes \psi_l \in
\widehat{D(H)}$, where $\{ f_i\}$, $\{\xi^j\}$ are dual bases in
$H,\widehat H$, respectively. 
\end{remark}

\begin{proposition}
\label{D(A) is factorizable}
The Drinfeld double $D(H)$ is factorizable in the sense of 
Definition~\ref{factorizability}.
\end{proposition}
\begin{proof}
First, observe that for any pair of dual bases $\{ f_i\}$ and $\{\xi^j\}$
as above and all $g\in H$ and $\psi\in \widehat{H}$ the element
$$
Q_{g\otimes \psi} = \sum_{k,l}\,f_kgf_l\otimes \xi^k\psi S(\xi^l)\in
H\otimes_k {\widehat H}^{op}
$$
belongs to $\widehat{D(H)}$. Indeed,
\begin{eqnarray*}
\lefteqn{\la Q_{g\otimes \psi},\, \phi \otimes zh \ra =} \\
&=& \la zh\1,\,\phi\1 \ra \la S(h\3),\,\phi\3 \ra \la g,\, \phi\2 \ra
    \la h\2,\psi \ra \\
&=& \la h\1,\, (\phi \actl z)\1 \ra \la S(h\3),\, (\phi \actl z)\3 \ra
    \la g,\, (\phi \actl z)\2 \ra \la h\2,\,\psi \ra \\
&=& \la Q_{g\otimes \psi},\, (\eps\actl z)\phi\otimes h\ra,
\end{eqnarray*}
for all $h\in H,\, \phi\in \widehat{H},\, z\in H_t$ and, likewise,
\begin{equation*}
\la Q_{g\otimes \psi},\, \phi \otimes yh \ra
= \la Q_{g\otimes \psi},\,(y\actr \eps) \phi\otimes h\ra.
\end{equation*}
Next, we compute
\begin{eqnarray*}
\lefteqn{(Q_{g\otimes \psi} \otimes \id)(\R_{21}\R)=}\\
&=& \sum_{i,j}\, Q_{g\otimes \psi} ([{\xi^j}\2 \otimes {f_i}\2])
    \,[\xi^i \otimes f_j] \, \la S({f_i}\1),\, {\xi^j}\1 \ra
    \la {f_i}\3, {\xi^j}\3  \ra \\
&=& \sum_{i,j}\, Q_{g\otimes \psi} ([{\xi^j}\2 \otimes f_i])\,
    [ S({\xi^j}\1) \xi^i {\xi^j}\3 \otimes f_j]  \\
&=& \sum_{i,j,k,l}\, [S({\xi^j})\xi^k\psi S(\xi^l)\xi^i \otimes
    f_jf_k g f_lf_i] \\
&=&  [ \eps\1'\la 1\2',\,\eps\2'\ra \psi \eps\1\la 1\2,\,\eps\2\ra
     \otimes 1\1' g 1\1 ]  \\
&=& [\psi \eps\1\la 1\2,\,\eps\2\ra \otimes g 1\1 ],
\end{eqnarray*}
where we used the identities
\begin{eqnarray*}
\sum_j\, S({\xi^j}\1) \la g,\,{\xi^j}\2 \ra \otimes {\xi^j}\3 \otimes f_j
&=& \sum_{i,j}\, S(\xi^i) \otimes \xi^j  \otimes f_i g f_j, \\
\sum_{i,j}\,S(\xi^i)\xi^j  \otimes f_if_j 
&=& \eps\1 \otimes 1\1 \la 1\2,\,\eps\2 \ra,
\end{eqnarray*}
that follow from (\ref{R-dual bases1}), (\ref{R-dual bases2}),
(\ref{dual bases}) and from axioms (\ref{S epst}),(\ref{S epss}) of a 
quantum groupoid.
Therefore,
\begin{eqnarray*}
\lefteqn{\la 1\2,\, \eps\2  \ra Q_{g1\1\otimes \psi\eps\1}(\R_{21}\R) =}\\
&=& [\psi\eps\1 \eps'\1 \otimes g 1\1'1\1] \la 1\2,\,\eps\2' \ra 
    \la 1\2',\,\eps\2 \ra \\
&=& [\eps\1 \otimes 1\1][\psi \otimes g] S([\eps\2 \otimes 1\2])
    = \Ad_1^l([\psi\otimes g]),
\end{eqnarray*}
Thus, we conclude from Lemma~\ref{expectations} that the map
$$
\widehat{D(H)} \own x \mapsto (Q_x \otimes \id)(\R_{21}\R) \in C_{D(H)}(D(H)_s)
$$
is surjective, i.e., $D(H)$ is factorizable.
\end{proof}
\end{section}

%%%%%%%%%%%%%%%%%%%%%%%%%%%%%%%%%%%%%%%%%%%%%%%%%%%%%%%%%%%%%%%%%%%%%%
%%%%%%%%%%%%%%%%   Ribbon quantum groupoids  %%%%%%%%%%%%%%%%%%%%%%%%%
%%%%%%%%%%%%%%%%%%%%%%%%%%%%%%%%%%%%%%%%%%%%%%%%%%%%%%%%%%%%%%%%%%%%%%
\begin{section}{Ribbon quantum groupoids}
\begin{definition}
\label{Ribbon WHA}
A ribbon quantum groupoid is a quasitriangular quantum groupoid 
$H$ with an invertible
central element $\nu\in H$ such that
\begin{equation}
\Delta(\nu) = \R_{21}\R(\nu\otimes \nu) \quad\mbox{and}\quad S(\nu)=\nu.
\end{equation}
The element $\nu$ is called a {\em ribbon element} of $H$.
\end{definition}
For an object $V$ of $\Rep(H)$ we define  the twist $\theta_V :V\to V$
to be the multiplication by $\nu$ :
\begin{equation}
\theta_V(v) =\nu\cdot v, \quad v\in V.
\end{equation}
\begin{proposition}
\label{twist}
Let $(H, \R,\nu)$ be a ribbon quantum groupoid. The family of
homomorphisms $\{\theta_V\}_V$ is a twist in the braided monoidal
category $\Rep(H)$ compatible with duality.
Conversely, if $\theta_V(v) =\nu\cdot v$ with $\nu\in H$ is a twist 
in $\Rep(H)$,
then $\nu$  is a ribbon element of $H$.
\end{proposition}
\begin{proof}
Since $\nu$ is an invertible central element of $H$, the homomorphism
$\theta_V$ is an $H$-linear isomorphism. The twist identity 
$c_{W,V}c_{V,W}(\theta_V\otimes\theta_W)=\theta_{V\otimes W}$
follows from the properties of $\nu$:
$$
c_{W,V}c_{V,W}(\theta_V\otimes\theta_W)(x) =
\R_{21}\R(\nu\cdot x\I \otimes \nu\cdot x\II) = \Delta(\nu)\cdot x
= \theta_{V\otimes W}(x),
$$
for all $x =x\I\otimes x\II\in V\otimes W$.
Clearly, the identity $\R_{21}\R(\nu\otimes \nu) =\Delta(\nu)$ is equivalent
to the twist property. It remains to prove that
$$
(\theta_V \otimes \id_{V^*})b_V(z) = (\id_V \otimes \theta_{V^*})b_V(z),
$$
for all $z\in H_t$, i.e., that
$$
\sum_i\,\nu z\1\cdot \xi^i \otimes z\2\cdot f_i =
\sum_i\,z\1\cdot \xi^i \otimes \nu  z\2\cdot f_i,
$$
where $\sum_i\,\xi^i \otimes f_i$ is the canonical element
in $V^*\otimes V$. Evaluating the first factors of the above equality
on an arbitrary $v\in V$, we get the equivalent condition :
$$
\sum_i\,(\nu z\1\cdot \xi^i)(v) z\2\cdot f_i =
\sum_i\, (z\1\cdot \xi^i)(v) \nu z\2\cdot f_i,
$$
which reduces to $z\2S(\nu z\1)\cdot v = S(z\1)\nu z\2\cdot v$.
The latter easily follows from the centrality of $\nu=S(\nu)$ and
properties of $H_t$.
\end{proof}

\begin{proposition}
\label{ribbon category}
The category $\Rep(H)$ is a ribbon category if and only if  
$H$ is a ribbon quantum groupoid.
\end{proposition}
\begin{proof}
Follows from   Propositions~\ref{monoidal with duality}, \ref{braiding},
and \ref{twist}.
\end{proof}

For any endomorphism $f$ of an object $V$ of $\Rep(H)$, we define, following 
(\cite{T2}, I.1.5), its {\em quantum trace}
\begin{equation}
\tr_q(f)=d_V c_{V,V^*}(\theta_Vf \otimes \id_{V^*})b_V
\label{q1-trace}
\end{equation}
with values in $\End(H_t)$ and the {\em quantum dimension} of $V$
by $\dim_q(V)=\tr_q(\id_V)$. The next lemma gives an explicit computation 
of $\tr_q$ and $\dim_q$ via the usual trace of endomorphisms.

\begin{proposition}
\label{quantum trace}
Let $(H, \R,\nu)$ be a ribbon quantum groupoid, $f$ be an endomorphism
of an object $V$ in $\Rep(H)$. Then 
\begin{equation}
\tr_q(f)(z) = \Tr(S(1\1)u\nu f)z1\2,\qquad \ dim_q(V)(z)=\Tr(S(1\1)u\nu)z1\2,
\end{equation}
where $\Tr$ is the usual trace of endomorphisms, and $u\in H$ is the Drinfeld 
element.
\end{proposition}
\begin{proof}
Since the trace of an endomorphism $h\in \End_k(H_t)$, in terms
of the canonical element $\sum_i f_i\otimes \xi^i\in V\otimes_k V^*$, is 
$\Tr(h)=\sum_i\,\xi^i(h(f_i))$, 
the definition of $\tr_q$ gives:
\begin{eqnarray*}  
\tr_q(f)(z)
&=& d_V c_{V,V^*}(\theta_Vf \otimes \id_{V^*})b_V(z) \\
&=& d_V( \sum_i\,\R\II z\2\cdot \xi^i \otimes \R\I\nu z\1\cdot f(f_i)) \\
&=& \sum_i\,(\R\II z\2\cdot \xi^i)(1\1 \R\I\nu z\1\cdot f(f_i))1\2  \\
&=& \sum_i \xi^i( S(\R\II z\2)1\1 \R\I\nu z\1\cdot f(f_i))1\2  \\
&=& \Tr(S(1\1)u\nu f)z1\2,
\end{eqnarray*} 
where we used formulas (\ref{dV}) and (\ref{bV})
defining $b_V$ and $d_V$.
\end{proof}

\begin{corollary}
\label{quantum trace for connected WHA}
Let $k$ be algebraically closed.
If $H$-module $H_t$ is irreducible (which happens exactly
when $H_t \cap Z(H) = k$, i.e., when $H$ is connected (\cite{N}, 3.11,
\cite{BNSz}, 2.4), then $\tr_q(f)$ and $dim_q(V)$ are scalars:
\begin{equation}
\label{q-trace}
\tr_q(f) = (\dim H_t)^{-1}\ \Tr(u\nu f), \qquad
\dim_q(V) = (\dim H_t)^{-1} \Tr(u\nu).
\end{equation}
\end{corollary}
\begin{proof}
An endomorphism of an irreducible module is multiplication by
a scalar, therefore, we must have $\Tr(S(1\1)u\nu f)1\2=\tr_q(f)(1)
=\tr_q(f) 1$.
Applying the counit to both sides and using that $\eps(1) = \dim H_t$,
we get the result.
\end{proof}
\end{section}

%%%%%%%%%%%%%%%%%%%%%%%%%%%%%%%%%%%%%%%%%%%%%%%%%%%%%%%%%%%%%%%%%%%%%%%
%%%%%%%%%%%%%%%   Modular quantum groupoids  %%%%%%%%%%%%%%%%%%%%%%%%%%
%%%%%%%%%%%%%%%%%%%%%%%%%%%%%%%%%%%%%%%%%%%%%%%%%%%%%%%%%%%%%%%%%%%%%%%
\begin{section}{Towards modular categories}

Let us first recall some definitions needed in this section. 
Let ${\mathcal V}$ be a ribbon $Ab$-category over $k$, i.e.,
such that all $\Hom(V,W)$ are $k$-vector spaces (for all objects
$V,W\in {\mathcal V}$) and both operations $\circ$ and $\otimes$ are
$k$-bilinear.

An object $V\in {\mathcal V}$ is said to be {\it simple} if any 
endomorphism of $V$ is multiplication by an element of $k$.
We say that a family $\{V_i\}_{i\in I}$ of objects of ${\mathcal V}$ 
dominates  an object $V$ of ${\mathcal V}$ if there exists a finite set
$\{V_{i(r)}\}_r$ of objects of this family (possibly, with repetitions)
and a family of morphisms  $f_r:V_{i(r)} \to V,g_r:V\to V_{i(r)}$ such that 
$\id_V = \sum_r f_r g_r$.

A modular category (\cite{T2}, II.1.4) is a pair consisting of
a ribbon $Ab$-category ${\mathcal V}$ and a finite family $\{V_i\}_{i\in I}$
of simple objects of ${\mathcal V}$ satisfying four axioms:
\begin{enumerate}
\item[(i)] There exists $0\in I$ such that $V_0$ is the unit object.
\item[(ii)] For any $i\in I$, there exists $i^*\in I$ such that $V_{i^*}$ is 
isomorphic to $V^*_i$.
\item[(iii)] All objects of ${\mathcal V}$ are dominated by the
family $\{V\}_{i\in I}$.
\item[(iv)] The square matrix $S=\{S_{ij}\}_{i,j\in I} =
\{ \tr_q(c_{V_i,V_j}\circ c_{V_j,V_i}) \}_{i,j\in I}$ is invertible over $k$
(here $\tr_q$ is the quantum trace in a ribbon category defined by 
(\ref{q1-trace})).
\end{enumerate}

If a quantum groupoid $H$ is connected and semisimple over an 
algebraically closed field, modularity of $\Rep(H)$ is equivalent to 
$\Rep(H)$ being ribbon and such that the matrix $S=\{S_{ij}\}_{i,j\in I} =
\{ \tr_q(c_{V_i,V_j}\circ c_{V_j,V_i}) \}_{i,j\in I}$, where $I$ is the set of 
all (equivalent classes of) irreducible representations, is invertible.

\begin{remark}
\label{recall integrals}
Recall that $h\in H$ is a left (resp.\ right) {\em integral} if $xh=\eps_t(x)h$
(resp.\ $hx=h\eps_s(x)$) for all $x\in H$ (\cite{BNSz}, 3.24).
A {\em Haar integral} is a two-sided integral $h$ which is normalized,
i.e., $\eps_t(h)=\eps_s(h)= 1$. Existence of a Haar integral in a quantum
groupoid $H$ is equivalent to $H$ being semisimple and possessing an invertible
element $g$ such that $S^2(x) =gxg^{-1}$ for any $x\in H$
and $\chi(g^{-1})\neq 0$
for all irreducible characters $\chi$ of $H$ (\cite{BNSz}, 3.27). 
\end{remark}

The following lemma extends a result known for Hopf algebras (\cite{EG}, 1.1).
\begin{lemma}
\label{factorizable implies modular}
Let $H$ be a connected, ribbon, factorizable quantum groupoid 
over an algebraically closed field  $k$, and assume that $H$ has
a Haar integral.  Then $\Rep(H)$ is a modular category.
\end{lemma}
\begin{proof}
Note that $H$ is semisimple by Remark~\ref{recall integrals}.
We only need to prove the invertibility of the matrix formed by
\begin{eqnarray*}
S_{ij} 
&=& \tr_q(c_{V_i,V_j}\circ c_{V_j,V_i}) \\
&=& (\dim H_t)^{-1} \Tr( (u\nu)\circ c_{V_i,V_j}\circ c_{V_j,V_i}) \\
&=& (\dim H_t)^{-1}  (\chi_j \otimes \chi_i)( (u\nu\otimes u\nu)\R_{21}\R),
\end{eqnarray*}
where $V_i$ are as above, $I = \{1,\dots n\},\ \{ \chi_j\}$ is a basis 
in the space $C(H)$ of characters 
of $H$ (we used above the 
formula (\ref{q-trace}) for the quantum trace).

Observe that the linear map $F : \phi\mapsto (\phi\otimes \id)(\R_{21}\R)$ 
takes any element of the form $\phi = \chi\actl u\nu$ (i.e., $\phi(h) =
\chi(u\nu h)~\forall h\in H$), where $\chi\in C(H)$,
into $Z(H)$. Indeed, for any such $\phi$ and all $\psi\in \widehat H,\,h\in H$ 
we  have, using the fact that $u\in C_H(H_s)$ 
(this follows from Lemma~\ref {yz}),
the properties of $\eps_s$ and $\eps_s^{op}$, the relation $\Delta^{op}(h)\R 
=\Delta(h) \R~(h\in H)$, and the centrality of $\chi$ : 
\begin{eqnarray*}
\la F(\phi)h,\,\psi \ra
&=& \la u\nu \R\II {\R'}\I,\,\chi\ra \la \R\I {\R'}\II h,\, \psi \ra \\
&=& \la u\nu \R\II {\R'}\I \eps_s^{op}(h\1),\,\chi\ra 
    \la \R\I {\R'}\II h\2,\, \psi \ra \\
&=& \la u\nu \R\II h\2 {\R'}\I S^{-1}(h\1),\,\chi\ra
    \la \R\I h\3 {\R'}\II,\, \psi \ra \\
&=& \la u\nu h\2 \R\II {\R'}\I S^{-1}(h\1),\,\chi\ra
    \la h\3 \R\I {\R'}\II,\, \psi \ra \\
&=& \la u\nu S(h\1) h\2 \R\II {\R'}\I,\,\chi\ra    
    \la h\3 \R\I {\R'}\II,\, \psi \ra \\
&=& \la u\nu \R\II {\R'}\I,\,\chi\ra \la h \R\I {\R'}\II,\, \psi \ra \\
&=& \la h F(\phi),\,\psi \ra,
\end{eqnarray*}
therefore $F(\phi)\in Z(H)$. Since $H$ is factorizable, we know
from Corollary~\ref{restriction of F} that the restriction
$$
F : \{ \phi\in \widehat H \mid \phi = \phi\circ\Ad_1^r \} \to C_H(H_s)
$$
is a linear isomorphism. Since $\chi\actl u\nu$ belongs to the subspace
on the left hand side, we have a linear isomorphism between $C(H)\actl u\nu$ 
and $Z(H)$, hence, there exists an invertible matrix $T=(T_{ij})$ representing 
the map $F$ in the bases of $C(H)$ and $Z(H)$, i.e., such that
$F(\chi_j\actl u\nu) = \sum_i\, T_{ij}e_i$. Then
\begin{eqnarray*}
S_{ij} 
&=& (\dim H_t)^{-1} \chi_i(u\nu F(\chi_j\actl u\nu)) 
    = (\dim H_t)^{-1} \sum_k\, T_{kj} \chi_i(u\nu e_k) \\
&=& (\dim H_t)^{-1} (\dim V_i) \chi_i(u\nu) T_{ij}. 
\end{eqnarray*} 
Therefore, $S = DT$, where $D = \mbox{diag}\{ (\dim H_t)^{-1} (\dim V_i) 
\chi_i(u\nu) \}$. If $g$ is an element from Remark~\ref{recall integrals}
then $u^{-1}g$ is an invertible central element of $h$ and $\chi_i(u^{-1})
\neq 0$ for all $\chi_i$. By Corollary~\ref{uv is central} $uS(u)=c$
is invertible central, therefore 
$\chi_i(u)= \chi_i(c)\chi_i(S(u^{-1})))\neq 0$. 
Hence, $\chi_i(u\nu)\neq 0$ for all $i$ and $D$ is invertible.
\end{proof}

\begin{example}
An example of a modular category can be constructed from
{\em elementary} quantum groupoids classified in \cite{NV1}, 3.2.
A quantum groupoid $H$ is called {\em elementary} if $H\cong M_n(k)$. 
Then it is determined, up to an
isomorphism, by one of its counital subalgebras
$$
H_t \cong \oplus_\alpha\, M_{n_\alpha}(k),
$$
$n_\alpha,~\alpha =1\dots N$ are positive integers, 
$n=\sum_\alpha\,n_\alpha^2$.
From this one can see that
$$ 
D(H) = [\widehat H^{op}\otimes H] = [\widehat H_s \widehat H_t\otimes H]
[\eps \otimes H] = H,
$$
where for subsets $A\subset \widehat H^{op}, b\subset H$ we set 
$[A\otimes B]=\{[a\otimes b]\in D(H)\vert a\in A,b\in B\}$.
Hence $H$ is the Drinfeld double of itself. 
The $R$-matrix of $H$ is
$$
R =\sum_{i,j,k,l\alpha}\, \frac{1}{n_\alpha}\, 
E_{jl\alpha}^{ik\alpha} \otimes E_{ij\alpha}^{kl\alpha},
$$
where $\{E_{ij\alpha}^{kl\beta} \}_{i,j=1\dots n_\alpha}$ is a system
of matrix units in $H$. 
Both the Drinfeld and the ribbon elements are equal to 1.
Thus, all the conditions of
Lemma~\ref{factorizable implies modular} are satisfied, the category $\Rep(H)$
is modular with a unique irreducible object.
\end{example}

\end{section}
%%%%%%%%%%%%%%%%%%%%%%%%%%%%%%%%%%%%%%%%%%%%%%%%%%%%%%%%%%%%%%%%%%%%%%%%%%
%%%%%%%%%%%%%%%%%% $C^*$-quantum groupoids and modular categories %%%%%%%%
%%%%%%%%%%%%%%%%%%%%%%%%%%%%%%%%%%%%%%%%%%%%%%%%%%%%%%%%%%%%%%%%%%%%%%%%%%
\begin{section}
{$C^*$-quantum groupoids and unitary modular categories}
\begin{definition}
\label{semisimple and C*}
A {\em $*$-quantum groupoid}
is a quantum groupoid over a field $k$ with involution,  
whose underlying algebra $H$ is equipped with an antilinear involutive
algebra anti-homomorphism $*:H\to H$ such that  
$\Delta\circ *=(*\otimes *)\Delta$.
A $*$-quantum groupoid is said to be {\em $C^*$-quantum groupoid}, if
$k=\mathbb{C}$ and $H$ is a finite-dimensional $C^*$-algebra, i.e., $x^*x=0$
if and only if
$x=0,\ \forall x\in H$. 
\end{definition}

Definition \ref{semisimple and C*} together with the uniqueness of the unit, 
counit and  antipode imply that 
$$
1^*=1,\quad \eps(h^*)= \overline{\eps(h)},\quad (S\circ *)^2 =\id
$$ 
for all $h$ in a $*$-quantum groupoid $H$. It is also easy to check the
relations 
$$
\eps_t(h)^*=\eps_t(S(h)^*),\ \eps_t(h)^*=\eps_t(S(h)^*),
$$
therefore, $H_t$ and $H_s$ are $*$-subalgebras. The dual 
$\widehat H$ is also a $*$-quantum groupoid with respect to the $*$-operation
\begin{equation}
\la\phi^*,h\ra=\overline{\la\phi,S(h)^*\ra}\qquad \mbox{ for all } 
\phi\in \widehat H,\ h\in H.
\end{equation}
The square of the antipode of a $C^*$-quantum groupoid is an inner
automorphism, i.e., $S^2(h) = ghg^{-1}$ for some $g\in H$.
It is easy to see that there is a unique such $g$ 
satisfying the following conditions (\cite{BNSz}, 4.4):
\begin{enumerate}
\item[(i)] $\tr(\pi_\alpha(g^{-1}))=\tr (\pi_\alpha(g))\neq 0$ 
for all irreducible representations $\pi_\alpha$ of $H$ 
(here $\tr$ is the usual trace on a matrix algebra);
\item[(ii)] $S(g)=g^{-1}$, and
\item[(iii)]
$\Delta(g) =(g\otimes g)\Delta(1) = \Delta(1) (g\otimes g)$.
\end{enumerate}
This element $g$ is called
the {\em canonical group-like element} of $H$.
\begin{remark}
\label{C*-dual}
(i) Any $C^*$-quantum groupoid satisfies the conditions of
Remark~\ref{recall integrals}, so it always possess a Haar integral.  

(ii) If $H$ is a $C^*$-quantum  groupoid, then its dual $\widehat H$ is also a
$C^*$-quantum  groupoid (see \cite{BNSz}, 4.5, \cite{NV1}, 2.3.10).
\end{remark}
 
Groupoid algebras and their duals give examples of commutative
and cocommutative
$C^*$-quantum groupoids if the ground field $k=\mathbb{C}$ 
(in which case $g^*=g^{-1}$ for all $g\in G$).

One can check that for a quasitriangular $*$-quantum groupoid
$\bR=\R^*$.
\begin{proposition}  
\label{*-double}
If $H$ is a $C^*$-quantum groupoid, then $D(H)$ is a quasitriangular 
$C^*$-quantum groupoid. 
\end{proposition}
\begin{proof}
First let us show that $\widehat{D(H)}$, equipped with a natural involution
$$
\la X^*,\, \phi\otimes h \ra
=\sum_k\, \overline{\la g_k,\,(S\phi)^* \ra}
             \overline{\la S(h)^*,\,\psi_k \ra},
$$
where $X =\sum_k\, g_k\otimes \psi_k \in \widehat{D(H)},
h\in H,\,\phi\in \widehat{H}$, is a $C^*$-subalgebra of
the tensor product $C^*$-algebra $H\otimes \widehat{H}^{op}$.
For this it suffices to show that $X^*|_J=0$, i.e.,
$\la X^*,\, \phi\otimes zh \ra = 
\la X^*,\,(\eps\actl z)\phi \otimes h \ra$,
$\la X^*,\, \phi\otimes yh \ra =
\la X^*,\,(y\actr \eps)\phi \otimes h \ra$ for $z\in H_t, y\in H_s$.
For instance, one computes:
\begin{eqnarray*}
\la X^*,\, \phi\otimes zh \ra
&=& \sum_k\, \overline{\la g_k,\,(S\phi)^* \ra}
             \overline{\la S(z)^*S(h)^*,\,\psi_k \ra},\\
\la X^*,\,(\eps\actl z)\phi \otimes h \ra
&=& \sum_k\, \overline{\la g_k,\, S((\eps\actl z)\phi)^* \ra}
             \overline{\la S(h)^*,\, \psi_k \ra},
\end{eqnarray*}
for all $z\in H_t$. The right-hand
sides of the above equations are equal since
$$
\overline{\la z,\,\eps\1\ra} S(\eps\2)^*
= \la S(z)^*,\eps\2\ra \eps\1.
$$
Similarly one gets the other relation.

To prove that the comultiplication of $\widehat{D(H)}$
is a $*$-homomorphism we compute
\begin{eqnarray*}
\Delta(X)^*
&=& \sum_{i,j,k}\, ( {g_k}\2^* \otimes {\xi^j}^*{\psi_k}\1^* {\xi^i}^*) \otimes
    (f_j^* {g_k}\1^* S(f_i)^* \otimes {\psi_k}\2^* )\\
&=& \sum_{i,j,k}\, ( {g_k}\2^* \otimes \xi^i {\psi_k}\1^* \xi^j) \otimes
    ( S(f_i) {g_k}\1^* f_j \otimes {\psi_k}\2^* ) = \Delta(X^*),
\end{eqnarray*}
where we use that $\sum_j\, (\xi^j)^* \otimes S(f_j)^* =
\sum_j\, \xi^j \otimes f_j$ for every pair of dual bases.
Thus, $\widehat{D(H)}$ is a $C^*$-quantum groupoid and so is $D(H)$ (see 
Remark~\ref{C*-dual}). 
\end{proof}
%\begin{remark}
%For another idea of the proof of this proposition see (\cite{BSz1}, 3).
%\end{remark} 
In \cite{EG} it was shown that a quasitriangular semisimple Hopf algebra is
automatically ribbon with ribbon element $\nu=u^{-1}$, where $u$ is 
the Drinfeld
element. We are able to get a similar result for $C^*$-quantum 
groupoids.

\begin{proposition}
\label{*-Ribbon WHA}
A quasitriangular $C^*$-quantum groupoid $H$ is automatically ribbon 
with ribbon element $\nu=u^{-1}g=gu^{-1}$, where $u$ is the Drinfeld
element from Definition~\ref{Drinfeld element} and $g$ is the canonical
group-like element implementing $S^2$.
\end{proposition}
\begin{proof}
Since $u$ also implements $S^2$ (Proposition \ref{elements u and v}), 
$\nu=u^{-1}g$ is central, therefore $S(\nu)$ is also central.
Clearly, $u$ must commute with $g$. The same Proposition gives
$\Delta(u^{-1}) = \R_{21}\R(u^{-1}\otimes u^{-1})$,
which allows us to compute
$$
\Delta(\nu) = \Delta(u^{-1})(g\otimes g) = 
\R_{21}\R(u^{-1}g\otimes u^{-1}g) = \R_{21}\R(\nu \otimes \nu).
$$
Propositions~\ref{properties of R} and \ref{elements u and v}
and the trace property imply that
\begin{eqnarray*}
\tr(\pi_\alpha(u^{-1})) 
&=& \tr(\pi_\alpha( \R^{(2)}S^2(\R^{(1)}) )) \\ 
&=& \tr (\pi_\alpha(S^3(\R^{(1)}) S(\R^{(2)}) ) = \tr(\pi_\alpha(S(u^{-1}))).
\end{eqnarray*}
Since $u^{-1} = \nu g^{-1}$ and $\nu$ is central, the above relation 
means that 
$$
\tr(\pi_\alpha(\nu))\tr(\pi_\alpha(g^{-1})) =
\tr(\pi_\alpha(S(\nu)))\tr(\pi_\alpha(g)),
$$
and, therefore, $\tr(\pi_\alpha(\nu)) = \tr(\pi_\alpha(S(\nu)))$ for
any irreducible representation $\pi_\alpha$, which shows that
that $\nu = S(\nu)$.
\end{proof}
\begin{corollary}
\label{C*-trace}
For a connected ribbon $C^*$-quantum groupoid $H$ we have:
$$
\tr_q(f) = (\dim H_t)^{-1} \Tr_V(g\circ f),\qquad \dim_q(V) = 
(\dim H_t)^{-1} \Tr_V(g).
$$
for any $f\in \End(V)$, where $V$ is an $H$-module. 
\end{corollary}
To define the (unitary) representation category $\URep(H)$ of a $C^*$-quantum
groupoid
$H$ we consider {\it unitary} $H$-modules, i.e., $H$-modules $V$ 
equipped with a scalar product 
$$
(\cdot,\cdot): V\times V\to \mathbb{C} \qquad
\mbox{ such that } \qquad 
(h\cdot v,w)=(v,h^*\cdot w)\ \forall h\in H, v,w\in V.
$$
The notion of a morphism in this category remains the same as in $\Rep(H)$.
The monoidal product of $V,W\in \URep(H)$ is defined as follows. We
construct a tensor product $V\otimes_\mathbb{C} W$ of Hilbert spaces and remark that
the action of $\Delta(1)$ on this left $H$-module is an orthogonal projection.
The image of this projection is, by definition, the monoidal product of
$V,W$ in $\URep(H)$. Clearly, this definition is compatible with the monoidal
product of morphisms in $\Rep(H)$.
 
For any $V\in \URep(H)$, the dual space $V^*$ is naturally identified 
($v\to \overline v$) with the conjugate Hilbert space, and under this 
identification we have 
$h\cdot\overline v= \overline{S(h)^*\cdot v}\ (v\in V, \overline v\in V^*$). 
In this way $V^*$ becomes a unitary $H$-module with scalar product 
$(\overline v, \overline w)=(w,gv)$, where $g$ is the canonical group-like 
element of $H$. 

The unit object in $\URep(H)$ is $H_t$ equipped with scalar product 
$(z,t)_{H_t}= \eps(zt^*)$ (it is known \cite{BNSz}, \cite{NV1} that
the restriction of $\eps$ to $H_t$ is a non-degenerate positive form). One can
verify that the maps 
$l_V,r_V$ and their inverses are isometries. For example, let us show that
the adjoint map for $l_V$ is exactly its inverse. We have:
$$
(l_V(1\1\cdot z\otimes 1\2\cdot v),w)
=(z\cdot v,w)\ (\forall z\in H_t,v,w\in V).
$$
On the other hand:
\begin{eqnarray*}
(1\1\cdot z\otimes 1\2\cdot v,l_V^{-1}w)
&=& (zS(1\1) \otimes 1\2\cdot v,S(1\1) \otimes 1\2\cdot w) \\
&=& \eps(zS(1\1)S(1\1)^*)(1\2^*1\2\cdot v,w) = (z\cdot v,w).
\end{eqnarray*}
Proposition \ref{monoidal with duality} implies that 
$\URep(H)$ is a monoidal category with duality
(see also \cite{BNSz}, Section 3).

\begin{remark}
\label{*-R-matrix} 
a) One can check that for
a quasitriangular $*$-quantum groupoid
the braiding is an isometry in $\URep(H)$:
$c^{-1}_{V,W}=c^{*}_{V,W}$.

b) For a ribbon $C^*$-quantum groupoid $H$,
the twist is an isometry in $\URep(H)$. 
Indeed, the relation $\theta_{V}^*=\theta_{V}^{-1}$ 
is equivalent to the identity $S(u^{-1}) = u^*$, 
which follows from Proposition~\ref{properties of R} and 
Remark~\ref{*-R-matrix}a).
\end{remark} 

A {\em Hermitian ribbon category} over the field $k$ with involution is an
$Ab$-ribbon category over $k$ endowed with a 
conjugation of morphisms $f\mapsto\overline{f}$ satisfying natural conditions
(see \cite{T2}, II.5.2):
\begin{eqnarray}
\label{1st line}
\overline{\overline{f}} &=& f,\qquad \overline{f+g}=\overline{f}+
\overline{g},\qquad  \overline{cf}=\overline c\overline{f} \quad
(c\in k), \\
\label{2nd line}
\overline{f\otimes g} &=& \overline{f}\otimes \overline{g},
\quad \overline{f\circ g}=\overline{g}\circ \overline{f},\qquad  
\overline{c_{V,W}}=(c_{V,W})^{-1}, 
\overline{\theta_{V}}=\theta_{V}^{-1}, \\
\label{herm}
\overline{b}_V &=& d_V\circ c_{V,V^*}(\theta_V\otimes id_{V^*}), \qquad
\overline{d}_V = (id_{V^*}\otimes \theta^{-1}_V)c^{-1}_{V^*,V}\circ b_V.
\end{eqnarray}
A {\em unitary ribbon category} is a Hermitian ribbon category over the
field $\mathbb{C}$ such that for any morphism $f$ we have $\tr_q(f\overline{f})
\geq 0$.

In a natural way we have a {\em conjugation}
of morphisms in $\URep(H)$. Namely, for any morphism $f:V\to W$ we define 
$\overline f:W\to V$ as $\overline f (w)=\overline{f^*(\overline w)}$ for
any $w\in W$. Here $\overline w\in W^*, f^*:W^*\to V^*$ is the standard dual
of $f$ (see \cite{T2}, I.1.8) and $\overline{f^*(\overline w)}\in V$.

\begin{lemma}
\label{hermit-ribbon}
Given a quasitriangular $C^*$-quantum groupoid $H,\ \URep(H)$ is a unitary 
ribbon $Ab$-category with respect to the  
above conjugation of morphisms.
\end{lemma}
\begin{proof}
Relations (\ref{1st line}) are obvious, (\ref{2nd line})
follows from Remarks~\ref{*-R-matrix}.

Let us prove relations (\ref{herm}).
On the one hand, for all $v\in V, \phi\in V^*, z\in H_t$ we have, using the 
definitions of $d_V,c_{V,V^*},\theta_V$, Propositions~\ref{properties of R}, 
\ref{elements u and v} and the notation $\omega_{v,\phi}(L)=(Lv,\phi)$
for a linear operator $L$ and two vectors $v,\phi$ of a Hilbert space:
\begin{eqnarray*}
\lefteqn{(d_V\circ c_{V,V^*}
(\theta_V\otimes id_{V^*})(v\otimes\phi),z)_{H_t} =} \\
&=& (d_V\circ c_{V,V^*}(gu^{-1}v\otimes\phi),z)_{H_t} \\
&=& (d_V[S(\R^{*(2)})\phi\otimes \R^{(1)}gu^{-1}v],z)_{H_t}  \\
&=& \eps[(1_{(1)}\R^{(1)}gu^{-1}v,S(\R^{*(2)})\phi)1_{(2)}z^*] \\
&=& (\omega_{v,\phi}\otimes\eps)
     [(S(\R^{(2)})\otimes 1)\Delta(1)(\R^{(1)}gu^{-1}\otimes z^*)] \\
&=& \omega_{v,\phi}[S(\R^{(2)})S(z^*)\R^{(1)}gu^{-1}]  \\  
&=& \omega_{v,\phi}[S(\R^{(2)})\R^{(1)}z^*gu^{-1}] = (z^*gv,\phi).
\end{eqnarray*}
And, on the other hand, using the definition of $b_V$, we compute :
\begin{eqnarray*}
(\overline{b}_V(v\otimes\phi),z)_{H_t}
&=& (v\otimes\phi,\sum_i z_{(1)}f_i \otimes
    S(z_{(2)})^*\xi^i)_{V\otimes V^*} \\
&=& \sum_i(v,z_{(1)}f_i)(S(z_{(2)})^*\xi^i, g\phi) \\
&=& \sum_i(z^*_{(1)}v,f_i)(\xi^i,S(z_{(2)})g\phi) \\
&=& (z^*_{(1)}v,S(z_{(2)})g\phi)=(v,zG\phi)=(z^*gv,\phi),
\end{eqnarray*}
whence the first part of (\ref{herm}) follows. To establish the
second part, note that  for  all $v\in V, \phi\in V^*, z\in H_t$ we have, 
using the definitions of $b_V,c_{V,V^*},\theta_V$, 
Propositions \ref{properties of R}, \ref{elements u and v} 
and the properties of $\nu$:
\begin{eqnarray*}
\lefteqn{((\id_{V^*}\otimes \theta^{-1}_V)c^{-1}_{V^*,V}
\circ b_V(z),\phi\otimes v)_{V^*\otimes V} =} \\
&=& ((\id_{V^*}\otimes \theta^{-1}_V)c^{-1}_{V^*,V}\sum_i z_{(1)}f_i
    \otimes S(z_{(2)})^*\xi^i),\phi\otimes v)_{V^*\otimes V}  \\
&=& ((\id_{V^*}\otimes \theta^{-1}_V)\sum_i S(\R^{*(1)})^*S(z_{(2)})^*\xi^i
    \otimes \R^{*(2)}z_{(1)}f_i),\phi\otimes v)_{V^*\otimes V}  \\
&=& \sum_i\,(\phi,gS(\R^{*(1)})^*S(z_{(2)})^*\xi^i)  
    (\nu^{-1}\R^{*(2)}z_{(1)}f_i,v)  \\
&=& \sum_i\,(S(z_{(2)})S(\R^{*(1)})g\phi,\xi^i)
    (f_i,\nu z^*_{(1)}\R^{(2)}v)  \\
&=& (S(z_{(2)})S(\R^{*(1)})\phi,\nu z^*_{(1)}\R^{(2)}v)
    =(\R^{*(2)}zS(\R^{*(1)})g\phi, \nu v) \\
&=& ((S^{-1}(\R^{(1)})\R^{(2)})^*S(z)g\phi,\nu v)
    =(S^{-1}(u)^*S(z)g\phi,\nu v) = (\phi,S(z^*)v).
\end{eqnarray*}  
On the other hand, using the definition of $d_V$, we obtain:
\begin{eqnarray*}
(\overline{d}_V(z),\phi\otimes v)_{V^*\otimes V}
&=& (z,d_V(\phi\otimes v))_{H_t} = (z,(1_{(1)}v,\phi)1_{(2)})_{H_t}  \\
&=& \overline{\eps((1_{(1)}v,\phi)1_{(2)}z^*)}
    =  \overline{(\omega_{v,\phi}\otimes \eps)(\Delta(1)(1\otimes z^*))} \\
&=& \overline{(\omega_{v,\phi}(S(z^*))}=(\phi,S(z^*)v). 
\end{eqnarray*}
The condition $\tr_q(f\overline f)=\Tr(gff^*)\geq 0$ for any morphism $f$
follows from Remark~\ref{*-R-matrix}b) and from the positivity of $g$.
\end{proof}
The next proposition extends (\cite{EG}, 1.2). 
 
\begin{theorem}
\label{rep of C* is modular}
If $H$ is a connected $C^*$-quantum groupoid, then $\URep(D(H))$ is a unitary 
modular category.
\end{theorem}
\begin{proof}
The proof follows from Lemmas~\ref{hermit-ribbon}, 
\ref{factorizable implies modular} and
Propositions~\ref{D(A) is factorizable}, \ref{*-double}.
\end{proof}
\end{section}

%%%%%%%%%%%%%%%%%%%%%%%%%%%%%%%%%%%%%%%%%%%%%%%%%%%%%%%%%%%%%%%%%%%%%%%%%%
%%%%%%%%%%%%%%%%%%%%  Appendix %%%%%%%%%%%%%%%%%%%%%%%%%%%%%%%%%%%%%%%%%%%
%%%%%%%%%%%%%%%%%%%%%%%%%%%%%%%%%%%%%%%%%%%%%%%%%%%%%%%%%%%%%%%%%%%%%%%%%%

\begin{section}{Appendix}
Here we collected some results on ribbon and modular quantum groupoids which
extend the corresponding facts for Hopf algebras.
\noindent
\textbf{1.} There is a procedure analogous to (\cite{RT1}, 3.4), that extends
any quasitriangular quantum groupoid $(H, \R,\bR)$ 
to a ribbon quantum groupoid in a canonical way. 
For this we need
\begin{lemma}[cf. (\cite{RT1}, 3.3)]
\label{properties of nu}
(i) A ribbon element $\nu$ satisfies
$$
\eps_t(\nu) = \eps_s(\nu) = 1 \quad \mbox{and} \quad \nu^2 = (vu)^{-1},
$$
where $u$ and $v$ are the elements defined in
Proposition~\ref{elements u and v}.
\newline \hskip 0.5cm (ii) If $\nu_1$ and $\nu_2$ are two ribbon elements of 
$(H, \R)$,
then $\nu_2 = E\nu_1$, where $E\in H$ is an invertible central element
such that $E= S(E) =E^{-1}$, $\Delta(E) = \Delta(1)(E\otimes E)$ (i.e.,
$E$ is group-like), and $\eps_t(E) = \eps_s(E) = 1$.
\end{lemma}
\begin{proof}
(i) The definition of counit implies :
\begin{eqnarray*}
\nu
&=& (\id\otimes\eps)\Delta(\nu)
    = \nu \R\II {\R'}\I \eps(\R\I {\R'}\II\nu) \\
&=& \nu \R\II {\R'}\I \eps(\eps_s(\R\I) {\R'}\II\nu) 
     = \nu {\R'}\I \eps(\eps_s({\R'}\II) \nu)  \\
&=& \nu S(1\2)\eps(1\1 \nu) =\nu S(\eps_t(\nu)),
\end{eqnarray*}
hence $\eps_t(\nu)=1$. We used here the identity $\eps(hg) =\eps(\eps_s(h)g),~
h,g\in H$ and Lemma~\ref{properties of R}. Similarly, $\eps_s(\nu)=1$.
Using the antipode property, we compute
\begin{eqnarray*}
1 &=& \eps_t(\nu) = m(\id\otimes S)\Delta(\nu) \\
&=& \R\II {\R'}\I S({\R'}\II) S(\R\I) \nu^2 \\
&=& v S^2(\R\II) S(\R\I) \nu^2  = vu \nu^2.
\end{eqnarray*}  
\newline (ii) Set $E= \nu_1^{-1}\nu_2$. Then $E$ is central and invertible,
$S(E)=E$, and from part (i) we conclude that $E^2=1$. Next,
$$
\Delta(E) = \bR\bR_{21}(\nu_1^{-1} \otimes \nu_1^{-1})
\R_{21}\R (\nu_2\otimes \nu_2) = \Delta(1)(E\otimes E).
$$
Applying the counit to both sides of the last equality, we get
$E = E\eps_t(E) = E\eps_s(E)$, i.e., $\eps_t(E) = \eps_s(E) =1$.
\end{proof}
\begin{proposition}
Let $\tilde{H} = H +H\nu$ be a central extension of $H$,
consisting of formal linear combinations $h+g\nu$ with $h,g\in H$.
Then $(\tilde{H}, \R, \nu)$ is a ribbon quantum groupoid with
operations
\begin{eqnarray*}
(h+g\nu)(h'+g'\nu) &=& (hh'+gg'(vu)^{-1}) + (hg' + gh')\nu, \\
\Delta(h+g\nu) &=& \Delta(h) + \Delta(g) \R_{21}\R(\nu\otimes\nu), \\
\eps(h+g\nu) &=& \eps(h) +\eps(g),\\
S(h+g\nu) &=& S(h) +S(g)\nu.
\end{eqnarray*}
Note that $\tilde{H}$ contains $H =\{ h+0\nu \mid g\in H\}$ as a 
quantum subgroupoid.
\end{proposition}
\begin{proof}
One verifies that $\Delta$ is a homomorphism exactly as in \cite{RT1}.
The properties of $\R$ and $\nu$ follow directly from definitions.
For the counit axiom we have, using the properties of counital maps,
Proposition~\ref{properties of R}, and Lemma~\ref{properties of nu} :
\begin{eqnarray*}
(\eps\otimes\id)\Delta(h+g\nu)
&=& h + \eps(g\1 \nu \eps_t(\R\II {\R'}\I))  g\2 \nu \R\I {\R'}\II \\
&=& h + g\2 S(\eps_s(g\1 \nu)) = h+g\nu, \\
(\id\otimes\eps)\Delta(h+g\nu)
&=& h + \R\II {\R'}\I g\1 \nu \eps(\eps_s(\R\I {\R'}\II)  g\2 \nu) \\
&=& h + S(\eps_t(g\2 \nu)) g\1 \nu = h+g\nu.
\end{eqnarray*}
Axioms (\ref{eps m}) and (\ref{Delta 1}) of 
Definition~\ref{finite quantum groupoid} can be verified by a direct
computation.

Next, we observe that $\eps_t(h+g\nu) = \eps_t(h) + \eps_t(g)$
and $\eps_s(h+g\nu) = \eps_s(h) + \eps_s(g)$.  The antipode axiom 
follows from the identity $\nu^2 vu =1$ provided by 
Lemma~\ref{properties of nu}:
\begin{eqnarray*}
m(\id\otimes S)\Delta(h+g\nu)
&=& \eps_t(h) + g\1 \R\II {\R'}\I S({\R'}\II) S(\R\I) S(g\2)\nu^2 \\
&=& \eps_t(h) + g\1 vu S(g\2)\nu^2 = \eps_t(h+g\nu),\\
m(S\otimes \id)\Delta(h+g\nu)
&=& \eps_s(h) + S(g\1) S(\R\I)S({\R'}\II) {\R'}\I \R\II  g\2 \nu^2 \\
&=& \eps_s(h) + S(g\1) vu g\2 \nu^2  = \eps_s(h+g\nu).
\end{eqnarray*}
The anti-multiplicative properties of the antipode follow from the facts
that $S(uv) = uv$ and $S(\nu) =\nu$. 
\end{proof}
\noindent
\textbf{2.} Let us establish a relation between modular quantum groupoids
and modular categories.
A morphism $f:V\to W$ in a ribbon $Ab$-category 
${\mathcal V}$ is said to be {\it negligible} if for any morphism $g:W\to V$ we 
have  $\tr(fg)=0$. ${\mathcal V}$ is said to be {\it pure} if all negligible 
morphisms in  this category are equal to zero. 
A purification procedure transforming any 
ribbon $Ab$-category into a pure ribbon $Ab$-category is described in 
(\cite{T2}, XI.4.2); this procedure transforms hermitian ribbon $Ab$-categories
into hermitian pure ribbon $Ab$-categories (\cite{T2}, XI.4.3).
We say that a family $\{V\}_{i\in I}$ of objects of ${\mathcal V}$ 
quasidominates  an object $V$ of ${\mathcal V}$ if there exists a finite set 
$\{V_{i(r)}\}_r$  of objects of this family (possibly, with repetitions) 
and a family of morphisms  $f_r:V_{i(r)} \to V,g_r:V\to V_{i(r)}$ such that 
$\id_V - \sum_r f_r g_r$ is negligible. If  ${\mathcal V}$ is pure, 
then quasidomination coincides with domination.
Let $(H,\R,\nu)$ be a ribbon quantum groupoid. Then an $H$-module $V$ of 
finite $k$-rank is said to be {\it negligible} if $\tr_q(f)=0$ for any
$f\in\End(V)$.  If $k$ is algebraically closed, then any 
irreducible $H$-module is a simple object of $\Rep(H)$.

\begin{definition}
\label{modular WHA}
A modular quantum groupoid consists of a ribbon quantum groupoid
$(H,\R,\nu)$ together with a finite family of simple $H$-modules 
of finite rank  $\{V\}_{i\in I}$ such that:
\begin{enumerate}
\item[(i)] for some $0\in I$, we have $V_0=H_t$, the unit object of $\Rep(H)$;
\item[(ii)] for each $i\in I$, there exists $i^*\in I$ such that $V_{i^*}$ is 
isomorphic to $V^*_i$;
\item[(iii)] for any $k,l\in I$, the tensor product $V_k\otimes V_l$ splits 
as a  finite direct some of certain $\{V\}_{i\in I}$ 
(possibly with multiplicities) and a negligible $H$-module; 
\end{enumerate}
To formulate the last condition, 
let $S_{i,j} =  \tr_q(c_{V_i,V_j}\circ c_{V_j,V_i}),~i,j\in I$, 
where the braiding $c_{V_i,V_j}$ was defined in \ref{QT WHA} 
and  the quantum trace $\tr_q$ in \ref{quantum trace}.
\begin{enumerate}
\item[(iv)] 
The square matrix $[S_{i,j}]_{i,j\in I}$ is invertible in $M_{|I|}(k)$.
\end{enumerate}
\end{definition}
For any modular quantum groupoid, we define a subcategory $\C$ of 
$\Rep(H)$ as follows. The objects of $\C$
are $H$-modules of finite rank quasidominated by 
$\{V\}_{i\in I}$ and morphisms are $H$-morphisms of such 
modules; all the operations in $\C$ are induced by the corresponding 
operations in $\Rep(H)$.
Now taking into account the results of the previous sections and repeating 
the proof of (\cite {T2}, XI.5.3.2), we have the first statement of 
the following

\begin{proposition}
\label{modular WHA-modular Rep}
If $(H,\R,\nu,\{V\}_{i\in I})$ is modular, then the 
subcategory $(\C, \{V\}_{i\in I})$ of $\Rep(H)$ is quasimodular in the 
sense of (\cite {T2}, XI.4.3).
Conversely, if $(\C, \{V\}_{i\in I})$ is
quasimodular, then $(H,\R,\nu,\{V\}_{i\in I})$ is modular.
\end{proposition}
The proof of the second statement follows directly from the comparison of 
\cite {T2}, XI.4.3 and the above definition of a modular quantum groupoid. 
Purifying $(\C, \{V\}_{i\in I})$ as in (\cite {T2}, XI.4.2), 
we get a modular category (\cite{T2}, II.1.4). 
\end{section}

%%%%%%%%%%%%%%%%%%%%%%%%%%%%%%%%%%%%%%%%%%%%%%%%%%%%%%%%%%%%%%%%%%%%%%%%%%
\bibliographystyle{amsalpha}

\end{document}